\documentclass[]{siamonline190516}
\usepackage[english]{babel}
\usepackage[utf8]{inputenc}
\usepackage{bm}
\usepackage{amssymb}
\usepackage{amsmath}
\usepackage{mathtools}

\usepackage{amsfonts}
\usepackage{graphicx}
\usepackage{epstopdf}

\graphicspath{{plots/}}
\ifpdf
\DeclareGraphicsExtensions{.eps,.pdf,.png,.jpg}
\else
\DeclareGraphicsExtensions{.eps}
\fi

\usepackage{algorithmic}
\usepackage{algorithm}


\newcommand{\yy}{\bm{y}}
\newcommand{\xx}{\bm{x}}
\newcommand{\dint}{\text{d}}
\newcommand{\cov}{\operatorname{cov}}
\newcommand{\divergence}{\operatorname{div}}
\DeclareMathOperator*{\esssup}{ess\,sup}
\DeclareMathOperator*{\essinf}{ess\,inf}

\newcommand{\cost}{\operatorname{cost}}

\title{A Dimension-adaptive Combination Technique for Uncertainty Quantification}
\author{Uta Seidler \thanks{Institute  for  Numerical  Simulation,  University  of  Bonn,  Friedrich-Hirzebruch-Allee  7,  53115  Bonn,  Germany  \email{griebel@ins.uni-bonn.de}, \email{seidler@ins.uni-bonn.de}} \and Michael Griebel \footnotemark[1] \thanks{Fraunhofer  Institute  for  Algorithms  and  Scientific  Computing  SCAI,  Schloss  Birlinghoven,  53754  SanktAugustin,  Germany}}

\begin{document}

\maketitle

\begin{abstract}We present an adaptive algorithm for the computation of quantities of interest involving the solution of a stochastic elliptic PDE where the diffusion coefficient is parametrized by means of a Karhunen-Lo\`eve expansion. The approximation of the equivalent parametric problem requires a restriction of the countably infinite-dimensional parameter space to a finite-dimensional parameter set, a spatial discretization and an approximation in the parametric variables. We consider a sparse grid approach between these approximation directions  in order to reduce the computational effort and propose a dimension-adaptive combination technique.
In addition, a sparse grid quadrature for the high-dimensional parametric approximation is employed and simultaneously balanced with the spatial and stochastic approximation. Our adaptive algorithm constructs a sparse grid approximation based on the benefit-cost ratio such that the regularity and thus the decay of the Karhunen-Lo\`eve coefficients is not required beforehand. The decay is detected and exploited as the algorithm adjusts to the anisotropy in the parametric variables. 
We include numerical examples for the Darcy problem with a lognormal permeability field, which illustrate a good performance of the algorithm: For sufficiently smooth random fields, we essentially recover the spatial order of convergence as asymptotic convergence rate with respect to the computational cost.
\end{abstract}


\maketitle


\section{Introduction}

In this article we present an adaptive algorithm for the computation of quantities of interest for problems in uncertainty quantification which is based on the sparse grid combination technique. To this end, as a model problem, we  focus  on the elliptic PDE 
\begin{align} \label{eq:modelproblem}
-\divergence( a(\xx,\omega) \nabla u(\xx, \omega))= f(\xx)   \quad \text{ in } D\subset \mathbb{R}^n, 
\end{align}
where uncertainty is introduced by the diffusion coefficient  $a(\xx, \omega)$, which is modeled by a random field and where equation \eqref{eq:modelproblem} holds for any realization $\omega$ in the sample space $\Omega$. The  spatial domain is denoted by ${D\subset\mathbb{R}^n}, {n=1,2,3}$, and the derivatives are with respect to the spatial variable $\xx$. 
Given information about the distribution of $a(\xx,\omega)$, our goal is to obtain statistical information of the solution to the PDE \eqref{eq:modelproblem}. In particular, we are interested in  $\mathbb{E}[\mathcal{F}(u)]$ where $\mathcal{F}$ either denotes a continuous linear functional $\mathcal{F}\colon H_0^1(D)\to \mathbb{R}$ or $\mathcal{F}(u)=u^p$ with $p=1,2$ for the computation of the first and second moment. 
Although this elliptic PDE  is a rather simple model, it is relevant in practice as it is at the core of many application problems ranging from thermal analysis in engineering over diffusion processes in biology and chemistry to subsurface and groundwater flow problems in hydrology. Moreover, for  water management and in environmental and energy problems (see, e.g. \cite{bear_modeling_2010}), where the permeability of the medium under consideration is usually unknown, it is important to understand the propagation of uncertainty through the flow model. 

Equation \eqref{eq:modelproblem} is often studied either with an affine (e.g. \cite{griebel_multilevel_2020, kuo_quasi-monte_2012, teckentrup_multilevel_2015}) or a lognormal diffusion coefficient (e.g. \cite{ernst_convergence_2018,kuo_multilevel_2017}).
In this paper,  we concentrate on the case of a  lognormal random field, where its logarithm is a Gaussian field characterized by a mean and a covariance function and which captures a large range of values. 
Moreover, we consider a parametrized version of \eqref{eq:modelproblem} 
such that the problem depends on a countable set of stochastic parameters $\yy = (y_1, y_2, y_3, \ldots)$ with 
$y_k \in \mathbb{R}$, for $k\in \mathbb{N}$. 
More precisely, we consider the logarithm of the diffusion coefficient represented in the Karhunen-Lo\`eve expansion (KL). 
Although this might not be always optimal as shown in \cite{bachmayr_sparse_2017}, it provides a straightforward parametrization with independent random variables.
When moments are known, the KL expansion is obtained by solving an eigenvalue problem. Otherwise, the field needs to be obtained from discretized samples. Either way, the  exact knowledge of the eigenvalues of the KL expansion  might not always be available.  

In order to compute the quantities of interest numerically, several approximation steps need to be applied. First, the  random diffusion coefficient $a$ has to be replaced by a finite-dimensional approximation.
For this, the KL expansion is truncated after $m$ terms. The truncated coefficient then depends on parameters $\yy_m=(y_1,\ldots, y_m)$ where  $m$ can be large. 
Furthermore, the expectation $\mathbb{E}[\mathcal{F}(u)]$ is replaced  by an $m$-dimensional integral
with a proper $m$-dimensional measure. This high-dimensional integration problem needs then to be approximated by a quadrature rule which 
relies on the evaluation of the integrand on a set of deterministic points or random samples $\yy_m^{(i)}$, $i=1,\ldots, N$. For each of these parameter values $\yy_m^{(i)}$, a deterministic PDE problem needs to be solved which requires a third discretization in the spatial variable. 
For the spatial approximation, we consider here a finite element method with piecewise linear elements since 
the solution to \eqref{eq:modelproblem} for fixed $\yy_m$ is  in $H^2(D)$ 
for sufficiently smooth diffusion coefficients and a right hand side in $L^2(D)$  (see \cite{charrier_finite_2013, graham_quasi-monte_2015}).

For the computation of the $m$-dimensional expectation,   Monte Carlo or Quasi Monte Carlo methods are often applied \cite{charrier_finite_2013, graham_quasi-monte_2015, harbrecht_multilevel_2016,teckentrup_further_2013} as their convergence rates with respect to the sample size are independent or nearly independent of the dimension of the stochastic parameter space. The convergence rates, however, are limited by a rate of $1/2$ and, depending on $m$, by a rate of a bit less than 1, respectively, even if the solution is smooth with respect to $\yy_m$. 
In fact, it has been shown in \cite{bachmayr_sparse_2017, harbrecht_multilevel_2016,hoang_n-term_2014} that in the lognormal case the dependence of the solution $u(\xx,\yy)$ on the parametric variables given by the KL expansion is analytic.
Furthermore, the influence of each stochastic parameter decreases due to the decay of the KL eigenvalues. 
Since these favorable regularity properties of the integrand cannot be exploited by conventional Monte Carlo or Quasi Monte Carlo quadratures, we instead choose for the parametric discretization a sparse grid type quadrature method based on univariate Gauss-Hermite quadrature rules 
such that  exponential convergence can indeed be obtained here. Other univariate quadrature rules  could also be considered, such as Genz-Keister or weighted Leja quadrature rules (see, e.g. \cite{chen_sparse_2018,narayan_adaptive_2014}). However, as we do not require nestedness of the univariate quadrature rules, we restrict our discussion for simplicity to the Gauss-Hermite quadrature which achieves the best convergence rate with respect to the number of quadrature points.
 Moreover,  we exploit the prevailing  anisotropy in the 
stochastic parameters  and apply a decay-obeying, dimension-adaptive sparse grid quadrature which has already been studied in several articles \cite{ babuska_stochastic_2010,haji-ali_novel_2018, nobile_anisotropic_2008, teckentrup_multilevel_2015}. This allows to harvest the favorable regularity properties of the integrand present in our specific problem.

Altogether, in this paper, we propose an adaptive algorithm that combines these different approximation steps, i.e., truncation, dimension-adaptive sparse grid quadrature, and finite element discretization in a multilevel manner. 
 We formulate a sparse grid approach with respect to these three discretization directions and, in this way, we balance the KL truncation, the  parametric and the spatial discretization with respect to the involved overall cost.
 Different aspects of the combination of these discretization steps have been considered separately before, but, at least to our knowledge, this is the first algorithm which adaptively steers {\em all} three discretization steps at once, i.e. the truncation, the quadrature, and the finite element method. A detailed discussion on the differences to other methods will be given in Section \ref{sec:dimadapt_UQ}. 
Similar to multilevel Monte Carlo methods and multilevel stochastic collocation methods \cite{barth_multi-level_2011, griebel_multilevel_2020, haji-ali_multi-index_2016, kuo_multilevel_2017, robbe_multi-index_2017, teckentrup_multilevel_2015,jakeman_adaptive_2020}, which involve different discretization levels for the spatial approximation and the quadrature, we can 
 benefit from  possible mixed regularity between the spatial and parametric variables.    
Note however that, in contrast to Monte Carlo and Quasi Monte Carlo methods, we have more than just one index to describe the parametric approximation as we use a sparse grid quadrature. Moreover, the number of discretization parameters for the quadrature is not fixed when we include the finite noise truncation as {\em additional} approximation direction in the multilevel approximation. In fact, the dimension of the   underlying quadrature problem varies since it depends on the truncation level. 
 To deal with this issue, we embed the parameters into an infinite-dimensional space and view the truncation as fixing higher parameters to the constant value $y_{m+j}=0$ for $j\geq 1$, which allows to easily modify the number of univariate quadrature rules. 
 
In order to balance simultaneously the error and cost of the finite noise truncation with each  univariate quadrature and with the spatial approximation, we use an adaptive variant of the well-known sparse grid combination technique \cite{Griebel.Schneider.Zenger:1992}. 
This has the advantage that we can employ straightforwardly standard finite element methods and quadratures. 
To this end, we modify the conventional dimension-adaptive combination technique \cite{gerstner_dimension_2003} to match our approximation steps.
This adaptive process is then steered by suitable indicators. We use here a reliable indicator based on the benefit-cost ratio for the adaptive construction of the sparse grid, see \cite{bungartz_sparse_2004}.  Therefore, in contrast to an a-priori construction, knowledge of the  regularity of the solution $u$ and of explicit decay rates of the KL eigenfunctions and eigenvalues is not necessary. 
In fact, it  suffices to be able to approximatively evaluate the KL truncation and to assume that the terms of the KL are sorted in decreasing order.
This  is similar to the {well-known} saturation assumption of adaptive finite element methods.
Consequently, our algorithm can be applied to different types of diffusion coefficients without any changes. 
It  properly balances the associated truncation error, the quadrature error and the finite element discretization error in a cost-versus-accuracy effective multilevel manner and adjusts itself to the anisotropy in the parametric variables.

The remaining part of this paper is structured as follows: We start by stating the considered problem and its parametrized version in Section \ref{sec:problemsetting} and continue with a separate description of each approximation step in Section \ref{sec:numericalmethods}. Subsequently, in Section \ref{sec:SG}, we focus on the general construction of sparse grids based on hierarchical surpluses and the combination technique, as well as the difference operators for our approximation. The dimension-adaptive sparse grid combination algorithm is then described in Section \ref{sec:dimadapt}. In Section \ref{sec:numex}, we present numerical examples which show the good convergence properties of our algorithm.


\section{Problem setting}\label{sec:problemsetting}

Let $(\Omega, \Sigma, \mathbb{P})$ be a probability space. 
As mentioned in the introduction, we want to compute the  quantities of interest 
\begin{align*}
	\mathbb{E}[\mathcal{F}(u)]= \int_\Omega \mathcal{F}(u(\xx, \omega))\dint \mathbb{P}(\omega)
\end{align*}
of the solution $u(\xx, \omega)$ to the random PDE 
\begin{align*}
	-\divergence(a(\xx,\omega)\nabla u(\xx, \omega))=f(\xx) \quad   \text{for } \xx\in D \quad \text{ for a.e.\ } \omega \in \Omega
\end{align*}
with $f\in L^2(D)$ and  zero Dirichlet boundary conditions for simplicity of notation. 
We assume that $a(\xx, \omega)$ is a lognormal random field, i.e., 
\begin{align} \label{eq:defalognormal}
a(\xx,\omega) = \exp(b(\xx, \omega))
\end{align}
where $b(\xx,\omega)$ is a Gaussian field expressed in form of a Karhunen-Lo\`eve expansion 
\begin{align} \label{eq:KLexpansion}
b(\xx,\omega) = \mathbb{E}[\log(a(\xx))]+ \sum_{k=1}^\infty \sqrt{\lambda_k} \psi_k(\xx)y_k(\omega)
\end{align}
with $\lambda_1 \geq \lambda_2 \geq \ldots$.
The functions $\{\psi_k\}$ are pairwise orthonormal functions in $L^2(D)$ and $y_k$ are independent standard normally distributed random variables, $y_k\sim \mathcal{N}(0,1)$, for which  we denote with $\rho_k$ the standard Gaussian density function.
In general, this expansion of $b$ can be obtained by computing the eigenpairs of the operator 
\begin{align*}
(C v) (\xx)=\int_D \cov_b(\xx,\bm{x'}) v(\bm{x'})\dint \bm{x'},
\end{align*}
whereby  $\cov_b(\xx,\bm{x'})$ is the continuous covariance function of the Gaussian field. Alternatively, it must be appropriately approximated from measurements. In this paper we do not require explicit knowledge of the KL eigenpairs and their decay rates, but we just need to be able to evaluate truncated and approximated versions of \eqref{eq:KLexpansion}.

Instead of the PDE \eqref{eq:modelproblem},  we consider a parametrized version in terms of the variables $y_k$ and set 
\begin{align*}
\yy \coloneqq (y_1, y_2, \ldots)  \in \mathbb{R}^{\mathbb{N}},
\end{align*} 
whose distribution is determined by the joint probability density function $\rho(\yy)=\prod_{k=1}^{\infty}\rho_k(y_k).$
The dependence of $a(\xx,\omega)$ on $\omega$  is substituted by a dependence on the parameters \linebreak  $\yy=(y_1,y_2,\ldots)$ and we write for the parametrized diffusion coefficient
\begin{align}\label{eq:a_parameterized}
a(\xx,\yy)=a_0(\xx) \exp\left(\sum_{k=1}^ \infty \sqrt{\lambda_k} \psi_k(\xx)y_k\right).
\end{align}
The parametrized version of the PDE is then for a fixed parameter $\yy\in\mathbb{R}^{\mathbb{N}}$ given by
\begin{align} \label{eq:ellitpicPDEparametrized}
		-\divergence(a(\xx,\yy)\nabla u(\xx,\yy))&=f(\xx) \qquad \xx \in D. 
\end{align}

For well-posedness of the problem, we need to ensure that the series in  \eqref{eq:a_parameterized} converges in $L^\infty(D)$ almost surely and that the moments are well-defined, i.e., $u\in L^2_\rho(\bm{\Gamma}, H^1_0(D))$. To this end, we assume that there exists a measurable set $\bm{\Gamma}\subset \mathbb{N}^\infty$ such that $\mu_{\rho}(\bm{\Gamma})=1$ and $a(\xx,\yy)$ converges in $L^\infty(D)$ for all $\yy\in\bm{\Gamma}$.
Furthermore, we define $a_{min}(\yy) = \essinf_{\xx \in D} a(\xx,\yy)$ and $\esssup_{\xx \in D} a(\xx,\yy)  = a_{max}(\yy)$ and assume that $a_{max}\in L^p_\rho(\bm{\Gamma})$ and $1/a_{min}\in L^p_\rho(\bm{\Gamma})$ for any $p\in (0,\infty)$. By the Lax-Milgram lemma, \eqref{eq:ellitpicPDEparametrized} admits for $\yy\in\bm{\Gamma}$ a unique solution $u(\yy)\in H_0^1(D)$  and satisfies 
\begin{align*}
\|u(\yy)\|_{H_0^1(D)} \leq \frac{C}{a_{min}(\yy)} \|f\|_{L^2(D)}.
\end{align*}
The expectation for the computation of the quantities of interest corresponds in the   \linebreak parametrized formulation then to the evaluation of the integral 
\begin{align} \label{eq:statparametrized}
\mathbb{E}\left[\mathcal{F}(u(\yy(\omega))) \right]=\int_\Omega \mathcal{F}(u( \yy(\omega))) \dint \mathbb{P}(\omega)=\mathbb{E}_{\yy}\left[\mathcal{F}(u(\yy))\right]= \int_{\bm{\Gamma}} \mathcal{F}(u( \yy)) \rho(\yy) \dint \yy 
\end{align}
which is well-defined for $\mathcal{F}(u)= u^p$, $p=1,2$ or for continuous linear functionals \linebreak  $\mathcal{F}: H_0^1(D) \to \mathbb{R}$. 
In the case where the KL eigenfunctions and -values are known explicitly, these assumptions can be ensured if the sequence 
\begin{align}\label{eq:defgamma}
\gamma_k=\sqrt{\lambda_k}\|\psi_k \|_{L^\infty(D)}
\end{align}
is summable, i.e., $\{\gamma_k\}_k\in \ell^1(\mathbb{N})$. 
In this case the set $\bm{\Gamma}$ is given by
\begin{align*}
\bm{\Gamma}=\{\yy \in \mathbb{R}^{\mathbb{N}} \colon \sum_{k=1}^\infty \gamma_k |y_k| < \infty \},
\end{align*}
which is $\mu_{\rho}$-measurable with $\mu_{\rho}(\bm{\Gamma})=1$ and $1/a_{min}, a_{max} \in L^p_\rho(\bm{\Gamma})$ for all $p\in (0,\infty)$ as shown in \cite{schwab_sparse_2011}.
However, note that also weaker assumptions are possible to ensure the well-posedness of the parametric problem (see, e.g. \cite{bachmayr_sparse_2017,charrier_strong_2012}) and we do not necessarily require summability of \eqref{eq:defgamma} in the following.

\section{Numerical approximations}\label{sec:numericalmethods}
This section focuses on the methods needed to numerically compute an approximation of \eqref{eq:statparametrized}.
In the following each approximation step  is addressed separately. We give a basic idea of the convergence behavior that can be expected if  sufficient regularity is provided, but refer to the literature for precise convergence rates and the regularity assumptions needed for  a precise convergence order. 
The combination of the different methods will be discussed in the later sections.

\subsection{Approximation by a finite-dimensional parameter set} 
First, since in practice one cannot work with the infinite representation \eqref{eq:KLexpansion} of the random diffusion field $a$, we approximate $a$ by its finite noise truncation, i.e., we cut off the series after $m$ terms. This way, we obtain
\begin{align*}
a_m(\xx,{\yy_m})=a_0(\xx)\exp\left(\sum_{k=1}^m \sqrt{\lambda_k} \psi_k(\xx) y_k\right),
\end{align*}
which depends on a finite but possibly high-dimensional parameter set. 
Furthermore, we denote with $u_m$ the solution of the PDE where the diffusion coefficient $a$ is replaced by $a_m$, i.e., for any ${\yy_m} \in \mathbb{R}^m$, the function $u_m({\yy_m})$ solves 
\begin{align*}
	-\divergence(a_m({\yy_m})\nabla u_m({\yy_m}))= f.
\end{align*}

For later purpose, we remark that the $m$-term truncation of the random diffusion coefficient $a$ is equivalent to replacing the dependence of $y_i$ for $i>m$ by the constant $0$. 
Thus, 
\begin{align}\label{eq:anchor}
	a_m(\xx,{\yy_m})=a(\xx,y_1, \ldots, y_m, 0,0, \ldots) \quad \text{ and } \quad 
	u_m(\xx,{\yy_m})=u(\xx,y_1, \ldots, y_m, 0,0, \ldots).
\end{align}

Note here that for $(y_1, \ldots,y_m, 0,0,\ldots )$ the series in \eqref{eq:a_parameterized} converges as it consists only of finitely many terms and the solution of the associated PDE is well-defined. 

Clearly, the error introduced by this truncation decreases with increasing $m$. The rate depends on the smoothness of the random field $a$ which is reflected in the decay rates of the eigenvalues and eigenfunctions of its Karhunen-Lo\`eve expansion.
 As shown in \cite{charrier_strong_2012, charrier_weak_2013, charrier_finite_2013, graham_quasi-monte_2015} also the solution $u_m$ converges to $u$ and error estimates for both the strong and weak error  are provided. These bounds depend on the decay properties of the coefficients in the KL expansion and require decay rates of $\lambda_k$ and $\psi_k$ in norms different to the $L^2(D)$-norm which are however in general difficult to obtain. 
But note here that, for our adaptive algorithm in Section \ref{sec:dimadapt}, we do not require these bounds to be available beforehand. 
Instead we will observe that our algorithm adjusts to different decay rates of the KL coefficients and determines the required truncation level automatically.
Hence, if the regularity of the random field $a$ is not known a-priori, we can still apply our algorithm in an a-posteriori fashion for a single specific $a$ in a similar way as adaptive finite element methods do.

\subsection{Finite element method}\label{sec:fem}
Next, we address the spatial discretization. 
For any fixed parameter ${\yy_m^{(i)}, i=1,\ldots,N}$, the solution $u_m(\cdot, {\yy_m^{(i)}})\in H_0^1(D)$ to the PDE \eqref{eq:ellitpicPDEparametrized} needs to be approximated. For this, we apply  a linear finite element method and subdivide the domain $D$ by a uniform and regular triangulation $\mathcal{T}_{h}$ with diameter $h$. 
Note here that we deliberately stick to uniform triangulations and do not consider spatial adaptivity since we are in this article merely interested in the interplay of truncation, plain discretization and quadrature.
For each $\yy\in \bm{\Gamma}$, we denote with the index $h$ the  finite element solution $u_{h}(\cdot, \yy)\in \mathcal{V}_h$ which  satisfies
\begin{align*}
\int_D a(\xx, \yy) \nabla u_{h}(\xx, \yy) \nabla v_h(\xx) \dint \xx= \int_D f(\xx) v_h(\xx) \dint \xx 
\end{align*}
for all $v_h\in \mathcal{V}_h=\{ v\in C(D) \colon v|_{\partial D}=0 \text{ and } v|_K\in \mathcal{P}_1 \; \forall K\in \mathcal{T}_h\}$,
where $\mathcal{P}_1$ is the space of linear functions. Moreover, we denote as $u_{h,m}(\cdot, \yy_m)\in \mathcal{V}_h$ the finite element solution of the truncated problem, i.e., 
\begin{align*}
\int_D a_m(\xx, \yy_m) \nabla u_{h,m}(\xx, \yy_m) \nabla v_h(\xx) \dint \xx= \int_D f(\xx) v_h(\xx) \dint \xx.
\end{align*}
For diffusion coefficients with sufficiently smooth realizations, the PDE solution satisfies  {$u(\cdot, \yy)\in H^2(D)$} a.s.\ and the  finite element approximation $u_h$ exhibits an error rate of the form 
\begin{align*}
	\| u(\xx, \yy)- u_h(\xx,\yy)\|_{H^1(D)} \leq  C  h \sqrt{\frac{a_{max}(\yy)}{a_{min}(\yy)}}  \|u(\xx,\yy)\|_{H^2(D)}.
\end{align*} 
For any $q<\infty$, this implies $\|u-u_h\|_{L^q_\rho(\bm{\Gamma},\,  H_0^1(D))} \leq C h $ and $\|u-u_h\|_{L^q_\rho(\bm{\Gamma},\,  L^2(D))} \leq C h^2 $  for the error in the spatial variable measured in the $L^2$-norm. These bounds also hold for the error $u_m-u_{h,m}$ with constants independent of $m$. 
For the computation of the second moment, the spatial norms are replaced by the $W^{1,1}(D)$-norm and the $L^1(D)$-norm, respectively. 
Hence, depending on the  quantity of interest and the  error norm, the error bound is of the order $\mathcal{O}(h)$ or $\mathcal{O}(h^2)$. 
For less regular diffusion coefficient without continuously differentiable realizations, the rate reduces to $\mathcal{O}(h^s)$ and $\mathcal{O}(h^{2s})$, respectively,  where $0<s\leq 1$ depends on the regularity of the random field. 
We refer to \cite{charrier_finite_2013, graham_quasi-monte_2015, teckentrup_further_2013} for more details on the spatial discretization error.

\subsection{Quadrature method}

As a third approximation step, we need to approximate the solution in the parametric variables and evaluate the integeral 
\begin{align}\label{eq:integrationwithdensity}
\bm{I}v =\int_{\mathbb{R}^m} v({\yy_m}) \rho({\yy_m}) \dint {\yy_m}
\end{align}
 where $v({\yy_m})=\mathcal{F}(u_{h,m}(\xx,{\yy_m}))$. 
Note that, with the truncation of the Karhunen-Lo\`eve expansion after $m$ terms, the integral is over the $m$-dimensional space $\mathbb{R}^m$. 
For this purpose, we consider quadrature rules of the type 
\begin{align*}
\bm{I} v\approx \bm{Q}_{N}v= \sum_{i=1}^{N} w_iv(\yy_m^{(i)})
\end{align*}
which are based on the evaluation of the integrand at $N$ quadrature points $\yy_m^{(i)}\in\mathbb{R}^m$  
that are taken into account with  weight $w_i$, $i=1, \ldots, N$. 
As $u_{h,m}$ is also a function in $\xx$, this requires a FEM-discretization on the mesh $\mathcal{T}_h$ and a solution of the resulting linear systems of equations for each of these $N$ parameter values. 

 We apply a Gaussian quadrature rule that can exploit the regularity of the integrand with respect to the parameters $y_1,y_2,\ldots$. For that, we combine $m$ one-dimensional Gauss-Hermite quadrature rules as the underlying density function for each parameter $y_k$ is the one-dimensional Gaussian. 
	For each parameter $y_k$ a univariate quadrature rule is considered which is based on the interpolation of the integrand among a set of deterministic points $\{y_k^{(i)}\}$  such that the highest degree of polynomial exactness is achieved. 
The quadrature points are chosen as the  $N^{(k)}$ roots of the Hermite polynomials of degree $N^{(k)}$, which are orthogonal with respect to the inner product
\begin{align*}
	(q,r)_{L^2_\rho(\mathbb{R})}= \int_\mathbb{R} q(y_k) r(y_k) \rho_k(y_k) \dint y_k,
\end{align*}
and the quadrature weights are obtained by integration of the interpolating polynomial, i.e., $w_i=(L_i, 1)_{L^2_{\rho_k}(\mathbb{R})}$ where $L_i$ denotes the $i$-th Lagrange polynomial with respect to the quadrature point set.

For the approximation of \eqref{eq:integrationwithdensity}, we combine these $m$ univariate Gauss-Hermite quadrature rules in a tensor product 
\begin{align*}
\bm{Q}_{(N^{(1)},N^{(2)},\dots, N^{(m)})}v=Q^{(1)}_{N^{(1)}}\otimes Q^{(2)}_{N^{(2)}}\otimes \ldots \otimes Q^{(m)}_{N^{(m)}} v 
\end{align*}
which requires the evaluation of $v$ at the  quadrature points  $(y_1^{(i_1)}, y_2^{(i_2)},\ldots, y_m^{(i_m)})$ for $i_k=1,\ldots N^{(k)}$ and $k=1,\ldots,m$. Here, $N=\prod_{k=1}^{m}N^{(k)}$.

If the integrand can be extended into a region in the complex plane, each univariate rule can achieve an exponential convergence rate. Hence, this choice of quadrature rule exploits the smooth dependence of the solution on the parameters $y_k$ for $k=1,\dots, m$, which is shown in  \cite{graham_quasi-monte_2015, harbrecht_multilevel_2016,hoang_n-term_2014}. 
However, the combination of the $m$ quadrature rules in a full isotropic  tensor product, i.e., with $N^{(1)}=N^{(2)}= \ldots= N^{(m)}$,   suffers from the curse of dimensionality as  the number of quadrature points grows exponentially in the dimension $m$. 
Instead, an anisotropic version, as discussed in \cite{babuska_stochastic_2010},  can be used which combines the $m$ rules with different approximation power  and  which can exploit the decreasing influence of the parameters $y_k$ on the PDE solution $u$ by choosing the numbers $N^{(k)}$ according to the anisotropy.
If the decrease of influence of the stochastic parameters $y_k$ is sufficiently fast, for example such that the errors of the univariate rules are balanced when ${N^{(k)}}= 2^{\alpha_k}$ with $\sum_{k=1}^\infty \alpha_k <\infty$, the curse of dimensionality can be broken. 
Furthermore, we combine such an  anisotropic approach with a sparse grid quadrature which leads to an anisotropic sparse grid quadrature similar to \cite{haji-ali_novel_2018, nobile_anisotropic_2008}. 
 The bounds in \cite{babuska_stochastic_2010, harbrecht_multilevel_2016,hoang_n-term_2014} on the derivatives of $u$ with respect to the parameters {$\yy_m$} justify that a sufficiently anisotropic sparse grid in the parametric variables is a sound approach. 
The construction of this anisotropic sparse grid quadrature is coupled with the KL-truncation and the spatial approximation which will become clearer in a moment.

\section{Sparse grids}\label{sec:SG}
The adaptive algorithm that we propose is based on a sparse grid approach which balances the three different approximation methods.
In the following, we  give a short overview of the general construction of sparse grids in $d$ dimensions with abstract operators before we address the adaptive combination technique for the computation of quantities of interest. 
For more details on sparse grids we refer to \cite{garcke_sparse_2012} and \cite{bungartz_sparse_2004}. The additional regularity between the spatial and parametric variables shown in \cite{ harbrecht_multilevel_2016, kuo_multilevel_2017} for sufficiently smooth random fields justifies this approach.

For each discretization direction we consider a sequence of numerical approximation operators $\{P_l\}_l$ with increasing approximation power such that the limit leads to the true solution. 

We define the difference operators between two approximations as
\begin{align*}
 \Delta_l \coloneqq \begin{cases}
 P_l -P_{l-1} &\quad \text{if } l\geq 1,\\
 P_0 &\quad \text{if } l=0,
 \end{cases}	
\end{align*}
such that $P_L$ can be reconstructed from the difference operators up to $L$, i.e.,
\begin{align*}
	P_L=\sum_{l=0}^L \Delta_l,
\end{align*}
which is just a telescoping sum.
For an approximation that involves $d$ different approximation directions, we combine the  $d$ operators  in a straightforward way and define the  anisotropic full grid approximation as
\begin{align*}
	\bm{P}_{\bm{L}} u= \bigotimes_{i=1}^d P^{(i)}_{L_i} u= P^{(1)}_{L_1}(P^{(2)}_{L_2}( \ldots P^{(d)}_{L_d}u ))= \sum_{\bm{l} \leq \bm{L}} \Delta_{\bm{l}}u 
\end{align*}
where $\bm{L}=(L_1,\ldots, L_d)$ indicates the approximation levels in each direction. The $d$-dimen-\linebreak sional difference operator is given by
\begin{align*}
\Delta_{\bm{l}}=\bigotimes_{i=1}^d \Delta_{l_{i}}^{(i)} \qquad \text{ for } \bm{l}=(l_1,\ldots, l_d)
\end{align*}  whereby the condition $\bm{l} \leq \bm{L}$ 
requires $l_i\leq L_i$ for all $i=1, \ldots, d$. 
Note at this point that the operators $P_{l_i}^{(i)}$ could in principle be of different type for different directions $i$,  i.e., it could be a spatial energy projection in one direction and quadrature rules in the other directions, as it will be done later on in subsection \ref{sec:SGUQ}.

A sparse grid approximation only includes differences for a subset $\mathcal{I}\subset \mathbb{N}^d$. For example, the regular sparse grid sums up the differences for $|\bm{l}|_1 \leq L+d-1$, which yields a good approximation if the solution has some additional smoothness in form of bounded mixed derivatives.
In the following, we apply a sparse grid with respect to a {\em general} index set $\mathcal{I}$ and define the sparse grid approximation by
\begin{align}\label{eq:spgeneralindex}
\bm{P}_{\mathcal{I}}^{\text{SG}} u= \sum_{\bm{l}\in \mathcal{I}}  \Delta_{\bm{l}} u.
\end{align}

The choice of the index set $\mathcal{I}$ plays here  a crucial role.
The optimal choice clearly depends on the contribution of each increment $\Delta_{\bm{l}} u $ and its associated cost. In \cite{bungartz_sparse_2004}, it has been shown that, for a prescribed cost, the selection of the index set can be formulated as a binary knapsack problem to which a solution is obtained by including the indices with the highest benefit-cost ratio, i.e., $\|\Delta_{\bm{l}} u\|$ in some norm  divided by the associated computational cost.
As a-priori estimates on the regularity of the solution and therefore on the size of $\Delta_{\bm{l}}u$ are not always available, we use the benefit-cost ratio  to construct the index set in a greedy way in our adaptive algorithm later on (see Section \ref{sec:dimadapt}).

It is known that the sparse grid approach allows to overcome the curse of dimensionality to some extent when the solution satisfies certain mixed smoothness properties.
A regular sparse grid, for example, reduces the cost compared to the full combination while nearly preserving the accuracy  such that the dependence of the error with respect to the computational cost on the dimension is only in the logarithmic term. 
In the case where the approximation power and the cost of the underlying methods are not equal, the approximation rate of the sparse grid using an anisotropic simplex for $\mathcal{I}$ is dominated by the slowest converging direction (see \cite{griebel_construction_2013}). Furthermore, in \cite{haji-ali_multi-index_2016}, an optimal index set has been identified for the particular case of exponential decay in the parametric approximation and of algebraic decrease in the spatial discretization error. In this case, the spatial discretization dominates the convergence. 
Note moreover that this aspect was also discussed in a more general context of approximation in UQ of the infinite dimensional problem in \cite{dung_hyperbolic_2016, dung__2018}.
Altogether, we can thus expect that our a sparse grid approach leads to a method where the direction with the weakest convergence rate dominates the overall convergence of the approximation.

\subsection{Combination technique}
For our algorithm we consider a somewhat different representation of the formula \eqref{eq:spgeneralindex}, the so-called combination technique, which  was  introduced in \cite{Griebel.Schneider.Zenger:1992}. While  efficient methods for the computation of the increments are necessary for \eqref{eq:spgeneralindex}, the combination technique relies on the combination of several anisotropic, full but low-order approximations.

The terms in \eqref{eq:spgeneralindex} are rearranged such that full grid combinations are recovered and the combination technique approximation for \eqref{eq:spgeneralindex} is defined as  
\begin{align} \label{eq:uct} 
\bm{P}^{\text{ ct}}_\mathcal{I} u=\sum_{ \bm{l}\in \mathcal{I}} \alpha_{\bm{l}} \bm{P}_{\bm{l}} u
\end{align} 
with  coefficients
\begin{align*}
\alpha_{\bm{l}}= \sum_{\bm{z}=\bm{0}}^{\bm{1}} \left(-1\right)^{|\bm{z}|_1} \chi^{\mathcal{I} }\left(\bm{l}+\bm{z}\right),
\end{align*}
where $\mathbf{1}=(1,\ldots, 1)$ and the characteristic function $\chi$ is given by
\begin{align*}
\chi^\mathcal{I}(\bm{l})=\begin{cases} 1 \quad \text{ if } \bm{l} \in \mathcal{I},\\
0 \quad \text{ otherwise.}\end{cases}
\end{align*}

In order to have an equivalence of the combination technique \eqref{eq:uct} and the sparse grid formulation based on the increments \eqref{eq:spgeneralindex}, we need to 
impose an admissibility condition on the index set $\mathcal{I}$ as in \cite{gerstner_dimension_2003}.
We require that the index set $\mathcal{I}$ is \textit{downward closed}, meaning it satisfies the following property:\\
\noindent
For any $\bm{l} \in \mathcal{I}$ it must hold that
\begin{align*}
\bm{l}-\bm{e}_j \in \mathcal{I} \quad \text{ for } 1\leq j \leq d \text{ with } l_j\geq 1
\end{align*}
with $\bm{e}_j$ being the $j$-th unit vector.\\
We will use the combination technique representation as it has  the advantage that it relies on anisotropic, full grid approximations and does not require the computation of the increments $\Delta_{\bm{l}}$. 
Hence, no  hierarchical structure needs to be assumed. Instead, the standard numerical methods discussed in Section \ref{sec:numericalmethods} can be directly applied, such as standard finite element solvers and standard quadratures, and no nestedness of the quadrature methods is needed. 
Additionally, the structure of \eqref{eq:uct} is intrinsically parallel as the approximations $\{\bm{P}_{\bm{l}}\}_{\bm{l}}$ can be computed independently.
Note however that  we will rely on the hierarchical operators $\Delta_{\bm{l}}$ for the adaptive construction of the sparse grid index set as explained later on. 

\subsection{Dimension-adaptive combination technique} \label{sec:dimadapt}
Next, we describe the adaptive construction of the index set $\mathcal{I}$. In the following, we state  the dimension-adaptive combination technique which was introduced for integration problems in \cite{gerstner_dimension_2003} and which will be modified to fit our problem setting later on in Section \ref{sec:dimadapt_UQ}. The basic algorithm with abstract operators $\bm{P}_{\bm{l}}$ is stated in Algorithm \ref{alg:generalalgo}.
The dimension-adaptive algorithm heuristically constructs a sequence of admissible index sets $\mathcal{I}^{(1)} \subset \mathcal{I}^{(2)} \subset \ldots \subset \mathcal{I}^{(t)}$ which indicate the increments that are included in the sparse grid formulation. 
Once such an index set has been found, the approximation is obtained using the  generalized combination technique formula \eqref{eq:uct}.

The algorithm starts off with one index in $\mathcal{I}^{(1)}$ and adds successively indices which do not destroy the admissibility of the index set.  To ensure this,  
two index sets are introduced $\mathcal{A}$ and $\mathcal{O}$ with $\mathcal{I}=\mathcal{A}\cup \mathcal{O}$ and $\mathcal{A}\cap \mathcal{O}=\emptyset$. The set $\mathcal{A}$ denotes the set of admissible neighbor  indices and the set $\mathcal{O}$ is referred to as the set of old indices. 
While the indices in $\mathcal{O}$ have already been chosen to be considered for the sparse grid, the index set $\mathcal{A}$ contains the indices in the neighborhood of $\mathcal{O}$, i.e., which result from a refinement along the axis, and which do not destroy the admissibility when added to $\mathcal{O}$. 

In each iteration of the dimension-adaptive combination technique, the index with the  highest contribution 
is selected among the indices in $\mathcal{A}$, removed from $\mathcal{A}$ and added to the set $\mathcal{O}$. This selection process is based on the benefit-cost ratio \cite{bungartz_sparse_2004} as mentioned before. Besides $\Delta_{\bm{l}}$, also the corresponding cost $c_{\bm{l}}$ is taken into account in the associated profit  indicator $\eta_{\bm{l}}$ which will be explained in detail later on. 
As a result of this process, the neighborhood for the index set $\mathcal{O}$ has changed such that  the  set $\mathcal{A}$ needs to be updated by adding the indices to $\mathcal{A}$ which still satisfy the admissibility condition.
In a final step, the full grid solutions are combined according to \eqref{eq:uct}. Since during the algorithm the solutions for $\bm{l}\in \mathcal{A}$ were computed for the profit indicator,  these indices are also included into the computation of the combination technique solution.

\begin{algorithm}
	\caption{Dimension-adaptive combination technique \cite{gerstner_dimension_2003}}
	\label{alg:generalalgo}
	\begin{flushleft}
		\textbf{Input:}  Tolerance $\varepsilon>0$\\
		\textbf{Output:} Index set $\mathcal{I}$ and approximation $\bm{P}^{\text{ ct}}_\mathcal{I} u$
	\end{flushleft}
	\begin{algorithmic}[1]
		\STATE $\bm{l}=\left(0,\ldots, 0\right)\in \mathbb{N}^d$
		\STATE $\mathcal{A}=\{\bm{l}\}$, $\mathcal{O}=\emptyset$
		\STATE compute local profit indicator $\eta_{\bm{l}}$ 
		\STATE $\eta=\eta_{\bm{l}}$
		\WHILE{ $\eta > \varepsilon$ }
		\STATE Select the index $\bm{j}\in \mathcal{A}$ with the largest profit $\eta_{\bm{j}}$ \label{algstep:selj}
		\STATE $\mathcal{A}=\mathcal{A}\setminus \{\bm{j}\}$, $\mathcal{O}= \mathcal{O}\cup \{\bm{j}\}$
		\STATE $\eta=\eta -\eta_{\bm{j}}$
		\FOR{$k=1, \ldots, d $}
		\STATE $\bm{l}= \bm{j}+ \bm{e}_k$
		\IF{$\bm{l}-\bm{e}_i \in \mathcal{O} \quad \forall 1 \leq i \leq d$ with $l_i>0$}
		\STATE $\mathcal{A}=\mathcal{A}\cup\{\bm{l}\}$
		\STATE Compute full grid approximation $\bm{P}_{\bm{l}} u$
		\STATE Compute local profit indicator $\eta_{\bm{l}}$
		\STATE $\eta=\eta + \eta_{\bm{l}}$
		\ENDIF
		\ENDFOR
		\ENDWHILE
		\STATE Combine approximations $\bm{P}_{\bm{l}} u$ for $\bm{l} \in \mathcal{I}=\mathcal{A} \cup \mathcal{O}$ according to \eqref{eq:uct} to obtain $\bm{P}^{\text{ ct}}_\mathcal{I} u$
		\RETURN $\mathcal{I}=\mathcal{A}\cup \mathcal{O}$ and $\bm{P}^{\text{ ct}}_\mathcal{I}u $
	\end{algorithmic}
\end{algorithm}


\section{Sparse Grids for Uncertainty Quantification}

In the previous section, we had a fixed dimension $d$ and described the basic so-called dimension-adaptive sparse grid approach for abstract approximation operators. Now, we generalize this approach to our uncertainty quantification problem and the computation of the quantities of interests \eqref{eq:statparametrized}. Here, the operators $P^{(i)}_{l_i}$ are of different type corresponding to the different approximation steps.  
The adaptation to uncertainty quantification  involves two main aspects: First, we employ a sparse grid for the quadrature of the parametric coordinates, but we also want to account for a possibly changing dimension due to the dynamically adapted truncation length. This will be done via the embedding of the parametric coordinates and corresponding indices into infinite-dimensional spaces, as already indicated in \eqref{eq:anchor}. 
Second, we include the spatial discretization in the sparse grid formulation by adding a spatial index to also incorporate the treatment of the associated deterministic PDEs. 
In the following, we first define the difference operators for the UQ setting. We then present the generalization of the dimension-adaptive combination technique to this problem setting and in the last subsection we discuss the UQ-specific profit indicators with which we steer the overall adaptive construction.

\subsection{Difference operators for UQ}\label{sec:SGUQ}
For our setting, we have three numerical approximations that need to be properly taken into account. We now define the difference operators for these three approximation directions.

First, we consider a sequence of truncation operators $\{T_m\}_m$ which map the PDE solution to the solution of the truncated problem, i.e., $T_m: u\mapsto u_m$, and we define the difference between two levels  by
\begin{align*}
\Delta^{(T)}_m u= \begin{cases}
u_{m}-u_{m-1} &\quad \text{if } m\geq 1,  \\
u_{0} &\quad \text{if } m=0.
\end{cases} 
\end{align*}

Second, for the spatial discretization, we employ a sequence of uniform grids on $D$ with width $h_{l_{\xx}}=\mathcal{O}(2^{-{l_{\xx}}})$ and consider  the operator that maps the function $u_m$ to the evaluation of $\mathcal{F}$ of its finite  element  approximation $\mathcal{F}(u_{h_{l_{\xx}},m})$. Due to the continuity of $\mathcal{F}$ we have  $\lim_{~{l_{\xx}}\to \infty} \mathcal{F}(u_{h_{l_{\xx}},m}) = \mathcal{F}(u_m)$.
Here, in the spatial direction, we have not just a projection onto the finite element subspace (which itself involves the FEM-discretization and the solution of the arising linear system of equations), but also the corresponding evaluation required for the specific quantity of interest. Only in the case of computing the first moment, i.e.,  $\mathcal{F}=\operatorname{Id}$, the operator corresponds to the projector onto the finite element space.
Consequently, the difference between two spatial discretizations is defined as
\begin{align*}
\Delta^{(\xx)}_{l_{\xx}} u_m= \begin{cases}
\mathcal{F}(u_{h_{l_{\xx}},m})-\mathcal{F}(u_{h_{{l_{\xx}}-1},m}) &\quad \text{if } {l_{\xx}}\geq 1, \\
\mathcal{F}( u_{h_0,m}) &\quad \text{if } {l_{\xx}}=0.\\
\end{cases}
\end{align*}
And finally, we approximate the expectation by a quadrature method for which we apply a sequence of quadrature rules that are themselves sparse grids. 
Here, given a truncation level {$m$}, we have {$m$} parametric variables. For each parameter $y_k$ for  $k=1,\ldots, {m}$, we consider 
a sequence of univariate Gauss-Hermite rules $\{Q^{(k)}_{l_k}\}_{l_k\geq0}$ with increasing number of quadrature points $\{N^{(y_k)}_{l_k}\}_{l_k\geq 0}$ and associated differences operators
\begin{align}\label{eq:quaddifference1d}
\Delta_{l_k}^{(Q_k)} = \begin{cases}
Q_{l_k}^{(k)}- Q_{l_k-1}^{(k)} &\quad \text{ if } l_k \geq 1,\\
Q_0^{(k)} &\quad \text{ if } l_k=0.
\end{cases}
\end{align} 
We combine the {$m$} one-dimensional difference operators and define the difference operator for the quadrature
\begin{align}\label{eq:simplesgQ}
\Delta_{\bm{l}}^{(Q)}=
\bigotimes_{k=1}^{{m}} \Delta_{l_k}^{(Q_k)}
\end{align}
with $\bm{l}$ being a {$m$}-dimensional array determining the level for each parametric variable.  
The summation of such differences over a set $\mathcal{I}_Q $ of indices $\bm{l}=(l_1, \ldots, l_m) \in \mathbb{N}^m$  then yields the sparse grid quadrature
\begin{align*}
\bm{Q}_{\mathcal{I}_Q}^{\text{SG}} v= \sum_{\bm{l}\in \mathcal{I}_Q}  \Delta_{\bm{l}}^{(Q)} v.
\end{align*}

In combination with the approximation of the random field $a(\xx,\yy)$ and the spatial discretization, the integration is however not with respect to a fixed high-dimensional parameter space, but should rather adapt itself properly, which gives rise to a varying number of parametric dimensions. To this end,  we exploit the fact that the $m$-term truncation of the KL expansion is equivalent to the constant evaluation at $y_k=0$ for $k>m$ (c.f.\ \eqref{eq:anchor}) {and embed the parameters and quadrature rules into an infinite dimensional setting}. 
A one-point Gauss-Hermite quadrature in $y_k$ for $k>m$ performs the same approximation and thus replaces the integrand by just a constant function in $y_k$. 
We therefore assume 
$N^{(y_k)}_0=1$,  by choosing $N^{(y_k)}_{l_k}=2^{l_k}$, and call the parameter dimension $k$ inactive if $l_k=0$. This lowest level of the univariate quadrature  then corresponds to a constant approximation and disregards the dependence of the random field on this $y_k$. Conversely, the parameter dimension $k$ is active when $l_k\geq 1$.

This allows us to represent the index indicating the quadrature levels by 
\begin{align}\label{eq:quadindexinfinity}
\bm{l}=(l_1, l_2, \ldots, l_m, 0, 0, \ldots ) \in \mathbb{N}^{\infty},
\end{align}
which is in principle independent of the parameter dimension, possesses  only finitely many non-zero entries, and altogether allows for  straightforward adaptivity in the parametric dimensions, i.e., in the truncation length.
The difference operator \eqref{eq:simplesgQ} is then equivalent to 
\begin{align*}
\Delta^{(Q)}_{\bm{l}} v=\Delta^{(Q)}_{(l_1,\ldots, l_m,0,0,\ldots)} v=\bigotimes_{k=1}^\infty \Delta_{l_k}^{(Q_k)} v =\bigotimes_{k=1}^m \Delta_{l_k}^{(Q_k)} v 
\end{align*}
with the univariate differences as in \eqref{eq:quaddifference1d}. The equality holds  because $Q_0^{(k)} v =v $ for \linebreak $k>m$.
Note also that the effective truncation level can be obtained from \eqref{eq:quadindexinfinity}  and is given by {$\max\{k\in \mathbb{N} \colon l_k\geq 1\}$}.

In order to also include the spatial discretization and to combine it with the parametric discretization, we extend the index  \eqref{eq:quadindexinfinity}  by one entry which refers to the spatial discretization and write
\begin{align}\label{eq:indexdef}
\tilde{\bm{l}}=(l_{\xx}, l_{y_1}, \ldots, l_{y_m}, 0,0,\ldots) \in \mathbb{N}^\infty.
\end{align}
Corresponding to this index, we set
\begin{align*}
	\bm{P}_{\tilde{\bm{l}}} u= Q_{l_{y_1}}^{(1)} \cdots Q_{l_{y_m}}^{(m)}\left(\mathcal{F}(u_{h_{l_{\xx}},m})\right),
\end{align*}
where $u_{h_{l_{\xx}},m}$ is the projection of $u_m$ onto the FEM-space with triangulation $\mathcal{T}_{h_{l_{\xx}}}$ (and involves thus discretization and solution of the resulting linear system of equations). Note here that the order of application of the quadrature, the application of $\mathcal F$, the FEM-projection (i.e., discretization and solution), and the truncation is not commutative and cannot be changed.

The generalized difference operator, which incorporates the parametric and spatial refinement as well as the truncation level, is then given by the combination of the three associated difference operators as
\begin{align}\label{eq:difference3d}
	\Delta_{\tilde{\bm{l}}} u &=\Delta_{\bm{l}}^{(Q)} \Delta^{(\xx)}_{l_{\xx}} \Delta^{(T)}_m u, 
\end{align}
where we use the notation $\Delta^{(\xx)}_{l_{\xx}} \Delta^{(T)}_mu$ for 
$\Delta^{(\xx)}_{l_{\xx}} \Delta^{(T)}_mu=\mathcal{F}(u_{h_{l_{\xx}}, m})-\mathcal{F}(u_{h_{l_{\xx}-1}, m})-\mathcal{F}(u_{h_{l_{\xx}}, m-1})+\mathcal{F}(u_{h_{l_{\xx}-1}, m-1}).
$ 
This difference operator \eqref{eq:difference3d} can further be simplified as follows: 
As the quadrature and truncation levels are not completely independent, we only need to consider indices where the effective truncation level coincides with the truncation level, i.e., \linebreak $m= \max\{k\in \mathbb{N} \colon l_{y_k}\geq 1\}$. In this case the differences can be written as 
\begin{align}\label{eq:naturaldifference}
\Delta_{\tilde{\bm{l}}} u &= \Delta_{(l_{y_1}, \ldots, l_{y_m}, 0, 0,\ldots)}^{(Q)} \Delta^{(\xx)}_{l_{\xx}} \Delta^{(T)}_m u	
= \bigotimes_{k=1}^\infty \Delta_{l_{y_k}}^{(Q_k)} \left(\Delta^{(\xx)}_{l_{\xx}} \Delta^{(T)}_m u\right) 		= 
		\bigotimes_{k=1}^m \Delta_{l_{y_k}}^{(Q_k)} \left(\Delta^{(\xx)}_{l_{\xx}} \Delta^{(T)}_m u\right)\nonumber  \\		
		&= \bigotimes_{k=1}^{m-1} \Delta_{l_{y_k}}^{(Q_k)} \otimes \Delta_{l_{y_m}}^{(Q_m)}\left(\mathcal{F}(u_{h_{l_{\xx}}, m})-\mathcal{F}(u_{h_{l_{\xx}-1}, m})-\mathcal{F}(u_{h_{l_{\xx}}, m-1})+\mathcal{F}(u_{h_{l_{\xx}-1}, m-1})\right).\nonumber
\end{align}
Note that the terms $u_{h_{l_{\xx}}, m-1}$ and $u_{h_{l_{\xx}-1}, m-1}$ are independent of the variable $y_m$. Because constant functions are integrated exactly by our quadrature rule, we obtain
\begin{align*}
&\ \Delta_{l_{y_m}}^{(Q_m)}\left(\mathcal{F}(u_{h_{l_{\xx}}, m-1})+\mathcal{F}(u_{h_{l_{\xx}-1}, m-1})\right)\\
=& \ Q^{(m)}_{l_{y_m}} \left(\mathcal{F}(u_{h_{l_{\xx}}, m-1})-\mathcal{F}(u_{h_{l_{\xx}-1}, m-1})\right) -Q^{(m)}_{l_{y_m}-1} \left(\mathcal{F}(u_{h_{l_{\xx}}, m-1})-\mathcal{F}(u_{h_{l_{\xx}-1}, m-1})\right)= 0
\end{align*}
since $l_{y_m}\geq 1$. 
Hence, we have for \eqref{eq:difference3d} the natural simplified representation
\begin{align}
\Delta_{\tilde{\bm{l}}} u
		&= \left(\bigotimes_{k=1}^{m-1} \Delta_{l_{y_k}}^{(Q_k)} \otimes \Delta_{l_{y_m}}^{(Q_m)}\right) \left(\mathcal{F}(u_{h_{l_{\xx}}, m})-\mathcal{F}(u_{h_{l_{\xx}-1}, m})\right)= \left(\bigotimes_{k=1}^m \Delta_{l_{y_k}}^{(Q_k)} \right) \Delta^{(\xx)}_{l_{\xx}} u_m,
	\end{align}
which only requires the quadrature differences of the $m$ active parametric variables $y_k$ (and the difference in the spatial variable $x$ of course).

\subsection{Dimension-adaptive algorithm for UQ} \label{sec:dimadapt_UQ}
We are finally in the position to modify the procedure described in Subsection \ref{sec:dimadapt} in order to match our approximation steps. 
To this end, we adapt the basic Algorithm \ref{alg:generalalgo} to now construct a sequence of  index sets  $\tilde{\mathcal{I}}^{(1)} \subset \tilde{\mathcal{I}}^{(2)} \subset \ldots \subset \tilde{\mathcal{I}}^{(t)}\subset \mathbb{N}^\infty$ which contain indices of the form \eqref{eq:indexdef}.
An index set $\tilde{\mathcal{I}}$ can then be used to obtain the sparse grid approximation 
\begin{align}\label{eq:SGUQ}
		\bm{P}_{\tilde{\mathcal{I}}}^{\text{SG}} u= \sum_{\tilde{\bm{l}}\in \tilde{\mathcal{I}}}  \Delta_{\tilde{\bm{l}}} u
		\qquad \text{ and } \qquad 
		\bm{P}^{\text{ ct}}_{\tilde{\mathcal{I}}} u=\sum_{ \tilde{\bm{l}}\in \tilde{\mathcal{I}}} \alpha_{\tilde{\bm{l}}} \bm{P}_{\tilde{\bm{l}}} u
\end{align}
analogously to \eqref{eq:spgeneralindex} and \eqref{eq:uct}.
	
The construction of $\tilde{\mathcal{I}}$ is also based on two sets $\tilde{\mathcal{A}}$ and $\tilde{\mathcal{O}}$, but we cannot store the infinite index arrays. However, since only a finite number of parametric dimensions are active at once, it suffices to store these active dimensions.
To facilitate the implementation, we store all parameter dimensions which have been activated in course of the algorithm even if some of them are inactive at the current iteration. The number of parameter dimensions that have been activated is denoted by $M$, which will get larger in the course of the algorithm. 
In addition to the indices corresponding to these $M$ activated quadratures dimensions, we store  the spatial discretization level and a fixed number  $\widehat{M}$ of  parameter dimensions that have not been activated yet to be able to increase the truncation level in future iterations of the algorithm.
We therefore store indices of the form
\begin{align*}
	\tilde{\bm{l}}=(l_{\xx}, l_{y_1},\ldots, l_{y_m}, 0, 0 \ldots, 0)\in \mathbb{N}^{1+M+\widehat{M}}
\end{align*}
with $m\leq M+\widehat{M}$ and use  a Boolean vector $act$ to indicate whether a variable has been activated.

Our new adaptive method for UQ problems, which is presented in Algorithm \ref{alg:dimadaptKTinfinite}, then proceeds similar to Algorithm \ref{alg:generalalgo}. 
It starts with one index in $\tilde{\mathcal{A}}$  and  adds successively indices with the highest profit indicator.  Here, the  set $\tilde{\mathcal{A}}$ includes indices which either refine the finite element discretization, refine the quadrature or raise the truncation level. 
The search for new  admissible indices is separated into two parts: 
Lines \ref{algstep:ex_forstart}-\ref{algstep:ex_forend} account for a
refinement of the finite element discretization  or the quadratures which have been activated
before, whereby including $\tilde{\bm{l}}+\tilde{\bm{e}}_1$ in $\tilde{\mathcal{A}}$ increases the level index $l_{\xx}$ and $\tilde{\bm{l}}+\tilde{\bm{e}}_k$ for $2\leq k \leq 1+M+\widehat{M}$ raises the level of the $k-1$-th direction, i.e., the $k$-th quadrature.
	
The lines \ref{algstep:ex_extendstart}-\ref{algstep:ex_extendend} become effective whenever $\tilde{\bm{l}}$ selected in line \ref{algstep:ex_selj} activates a new variable which corresponds to a higher truncation level. In this case, the Boolean vector needs to be adjusted and the  number of stored parametric level is raised by one. This ensures that always $\widehat{M}$ non-activated entries are included and the truncation level can be increased in subsequent iterations. Moreover, the index arrays  which were added in previous steps of the algorithm need to be extended by one entry which is set to zero. 
	
\begin{algorithm}
\caption{Adapted dimension-adaptive combination technique for UQ}
\label{alg:dimadaptKTinfinite}
\begin{flushleft}
\textbf{Input:}  Tolerance $\varepsilon>0$\\
\textbf{Output:} Index set $\tilde{\mathcal{I}}$ and approximation $\bm{P}^{\text{ct}}_{\tilde{\mathcal{I}}} u$
\end{flushleft}
    \begin{algorithmic}[1]
    \STATE $\tilde{\bm{l}}=\left(0,\ldots,  0\right)\in \mathbb{R}^{1+\widehat{M}}$, $M=0$, $act=\left(0, \ldots, 0\right)\in \{0,1\}^{\widehat{M}}$
    \STATE $\tilde{\mathcal{A}}=\{\tilde{\bm{l}}\}$, $\tilde{\mathcal{O}}=\emptyset$
    \STATE Compute local profit indicator $\eta_{\tilde{\bm{l}}}$ 
    \STATE $\eta=\eta_{\tilde{\bm{l}}}$
    \WHILE{ $\eta > \varepsilon$ }
        \STATE Select index $\tilde{\bm{l}}\in \tilde{\mathcal{A}}$ with largest profit $\eta_{\tilde{\bm{l}}}$ \label{algstep:ex_selj}
        \STATE $\tilde{\mathcal{A}}=\tilde{\mathcal{A}}\setminus \{\tilde{\bm{l}}\}$, $\tilde{\mathcal{O}}= \tilde{\mathcal{O}}\cup \{\tilde{\bm{l}}\}$
        \STATE $\eta=\eta -\eta_{\tilde{\bm{l}}}$
        \FOR{$k=1, \ldots, 1+M+\widehat{M} $} \label{algstep:ex_forstart}
        	\STATE $\tilde{\bm{j}}= \tilde{\bm{l}}+ \tilde{\bm{e}}_k$
        	\IF{$\tilde{\bm{j}}-\tilde{\bm{e}}_i \in \tilde{\mathcal{O}} \quad \forall 1 \leq i \leq 1+M+\widehat{M}$ with $\tilde{j}_i>0$}
        		\STATE $\tilde{\mathcal{A}}=\tilde{\mathcal{A}}\cup\{\tilde{\bm{j}}\}$
        		\STATE Compute the  approximation $\bm{P}_{\tilde{\bm{j}}} u$
        		\STATE Compute   local profit indicator $\eta_{\tilde{\bm{j}}}$
        		\STATE $\eta=\eta + \eta_{\tilde{\bm{j}}}$
        	\ENDIF
        \ENDFOR \label{algstep:ex_forend}
        \IF{ it exists $n\in \{1\ldots,M+\widehat{M}\}$ s.t. $l_{y_n}>0$ and $act_n=0$} \label{algstep:ex_extendstart}
        \STATE $act_n=1$, extend $act$ 
         \STATE $M=M+1$
        \STATE $\tilde{\bm{j}}=\tilde{\bm{e}}_{1+M+\widehat{M}}$ 
        \STATE $\tilde{\mathcal{A}}=\tilde{\mathcal{A}}\cup \{\tilde{\bm{j}}\}$ \label{algstep:ex_includel}
        \STATE Extend indices in $\tilde{\mathcal{A}}$, $\tilde{\mathcal{O}}$ by $0$
        \STATE Compute the  approximation $\bm{P}_{\tilde{\bm{j}}} u$
        \STATE Compute local profit indicator $\eta_{\tilde{\bm{j}}}$
        \STATE $\eta=\eta+\eta_{\tilde{\bm{j}}}$
        \ENDIF \label{algstep:ex_extendend}
    \ENDWHILE
    \STATE Combine approximations $\bm{P}_{\tilde{\bm{l}}} u$ for $\tilde{\bm{l}} \in \tilde{\mathcal{A}} \cup \tilde{\mathcal{O}}$ according to \eqref{eq:SGUQ} to obtain $\bm{P}_{\tilde{\mathcal{I}}}^{\text{ ct}}u$ 
        \RETURN $\tilde{\mathcal{I}}=\tilde{\mathcal{A}}\cup \tilde{\mathcal{O}}$ and approximation $\bm{P}^{\text{ ct}}_{\tilde{\mathcal{I}}} u $
    \end{algorithmic}
\end{algorithm}
		
Note here that it would be sufficient to consider $\widehat{M}=1$. By setting $\widehat{M}>1$, the algorithm does not only consider raising the truncation level by one, but can consider the $\widehat{M}$ forward neighbors as well. This avoids taking into account a too small number of stochastic variables as a parametric variable can be activated even if not all parameters up to this parameter have been activated before. In our experiments it turned out that this approach with value $\widehat{M}=5$ was more robust than with the choice $\widehat{M}=1$. 
This lack of downward closedness in the truncation direction does not pose any problems for the combination technique due to the structure of the differences, see also \eqref{eq:naturaldifference}.

Let us comment here  on the effect of choosing a different sequence of KL truncation levels: Instead of increasing the truncation level by one at a time, a doubling of the number of terms in the KL would rather resemble the geometric refinement of the spatial discretization. This would be straightforward for Monte Carlo quadratures. However, the truncation level and a sparse grid quadrature are not completely independent: Unlike Monte Carlo methods the sparse grid quadrature relies on univariate quadrature rules and therefore the dimension intrinsically divides into steps of one. Thus, we here propose an increment by one. Nevertheless, doubling the truncation could be realized together with a sparse grid quadrature as follows: Let $m$ indicate the truncation level with $2^m$ terms in the KL expansion. Then raising the level from $m-1$ to $m$ includes $2^{m-1}$ new parameter variables $y_k$ which need to be considered separately by a univariate rule. One way to do this is to employ a quadrature with quadrature level $l_{y_k}=1$ for all new variables $y_k$, $k=2^{m-1}+1,\ldots, 2^m$, which activates all new variables at once. But the cost of evaluating the quadrature grows with $\mathcal{O}(2^m)$ which inhibits including high dimensions.
Instead we allow some, but not all of the $y_k$, $k=2^{m-1}+1,\ldots, 2^m$ to be inactive. Raising the truncation level would then imply to include the new forward neighbors, $\tilde{\bm{e}}_k$ for $k=2^{m-1}+2, \ldots, 2^m+1$, in the set $\tilde{\mathcal{A}}$ in the algorithm and the length of index sets to be extended by $2^{m-1}$  instead of one. Note that any quadrature with $l_{y_k}=0$ for all $k>M$ for some $M\in \{2^{m-1}+1, \ldots, 2^m\}$ has the effective truncation level $M$. Hence in practice, any truncation level is implicitly assumed by this choice, and doubling the included KL terms corresponds to increasing the buffer variables $\hat{M}$ depending on $m$. Because of the anisotropy in the parametric variables, the outcome is nearly the same for both strategies of truncation refinement and there is no benefit by doubling the truncation level.

Next, let us remark on the differences to other adaptive algorithm for random elliptic PDEs based on \cite{gerstner_dimension_2003} which consider certain aspects of the approximation separately.  
An adaptive sparse grid quadrature method for parametric PDEs  with fixed truncation level has been for example studied in \cite{beck_optimal_2012, nobile_convergence_2016}. In \cite{garcke_adaptive_2016,chen_sparse_2018,ernst_convergence_2018, chkifa_high-dimensional_2014}, 
 adaptive quadrature and interpolation algorithms based on sparse polynomial approximation have been considered in an infinite dimensional parameter space for the lognormal as well as for the affine case. In these approaches, like in our algorithm, the number of active variables is increased throughout the process and also the concept of buffer variables was introduced. However, we additionally include the spatial discretization in our approximation. This allows us to improve the overall computational complexity and to benefit from mixed smoothness between the spatial and parametric variables. Furthermore, in our algorithm, the extension of the index set by $\hat{M}$ terms is interpreted as possible refinement of the KL truncation which better balances the cost and error of the truncation with the other discretizations. This perspective is also included in the choice of our profit indicator. 
In the context of multilevel methods, adaptive algorithms based on the benefit-cost ratio have been proposed as well \cite{haji-ali_multi-index_2016, jakeman_adaptive_2020, robbe_multi-index_2017}. But these methods only focus on the interplay of the  spatial/deterministic discretization and parametric approximation and require a fixed number of discretization parameters. We instead additionally include the truncation level  which leads to a dependence of two discretization parameters as the truncation also influences the parametric discretization. 
Note that an adaptive algorithm which adjusts the truncation level in a multilevel Monte Carlo method has been proposed in \cite{owen_dimension-adaptive_2018}, but does not encounter this problem as it uses the plain Monte Carlo quadrature. We use a different approach based on Gaussian quadrature which leads to a more complicated algorithm but allows to obtain convergence rates much better than $1/2$. 
Although certain aspects of our algorithm may not be particularly novel, the combination of the spatial discretization, KL truncation and sparse grid quadrature with varying number of dimensions in a sparse grid approximation surely is. 
This perspective is also reflected in the  choice of our specific  indicator which will be discussed in the next section.

\subsection{Profit indicator}
\label{sec:profit}

The crucial part of the adaptive algorithm is the refinement strategy which selects the index to be added to $\tilde{\mathcal{O}}$. This selection process is controlled by a local profit indicator  $\eta_{\tilde{\bm{l}}}$ which should estimate the contribution of including the increment $\Delta_{\tilde{\bm{l}}}u$ in the sparse grid combination formula \eqref{eq:SGUQ}. 

Since the global error can be bounded by
\begin{align*}
\|\mathbb{E}\left[\mathcal{F}(u)\right]-\bm{P}_{\tilde{\mathcal{I}}}^{\text{ct}} u\|_{\mathcal{X}} \leq \sum_{\tilde{\bm{l}}\notin \tilde{\mathcal{I}}} \| \Delta_{\tilde{\bm{l}}} u \|_{\mathcal{X}},
\end{align*}
 we use 
\begin{align*}
E_{\tilde{\bm{l}}} = \left \| \Delta_{\tilde{\bm{l}}} u \right \|_\mathcal{X} 
\end{align*}
as a local error estimator for adding $\tilde{\bm{l}}$ to $\tilde{\mathcal{I}}$.
The increments are as discussed in Section \ref{sec:SGUQ} and the norm $\|\cdot \|_\mathcal{X}$ depends on the space in which the error is measured. 
Note at this point that, due to the upper bound by the triangle inequality, the local error indicator $E_{\tilde{\bm{l}}}$ is {\em reliable}. It is however not necessarily {\em efficient} in the terminology of error estimators for adaptive finite element methods (see, e.g. \cite{ainsworth_posteriori_2000}) as we have no lower bound $ c \sum_{\tilde{\bm{l}} \notin \tilde{\mathcal{I}}} \|\Delta_{\tilde{\bm{l}}} u\|_\mathcal{X} \leq  \|\mathbb{E}\left[\mathcal{F}(u)\right]-\bm{P}_{\tilde{\mathcal{I}}}^{\text{ct}} u\|_{\mathcal{X}}$, $0 <c\leq 1$.

Note here that the evaluation of  $E_{\tilde{\bm{l}}}$ requires several anisotropic grid solutions. However, we only need to compute $\bm{P}_{\tilde{\bm{l}}}u$ and can reuse the other solutions that have been previously computed as the index set needs to be downward closed in the quadrature and finite element levels. 
Moreover, taking into account the computational cost, the local error estimator $E_{\tilde{\bm{l}}}$ might not be a well-suited  refinement indicator, as including other indices might lead to a similar error reduction while having much lower computational cost. 
Therefore, we follow  \cite{bungartz_sparse_2004, gerstner_dimension_2003} and use instead  the benefit-cost ratio 
\begin{align*}
\eta_{\tilde{\bm{l}}}= \frac{E_{\tilde{\bm{l}}}}{c_{\tilde{\bm{l}}}}
\end{align*}
as refinement indicator 
where $c_{\tilde{\bm{l}}}$ denotes the cost of the evaluation of $\bm{P}_{\tilde{\bm{l}}}u$. As cost model we use the product of the cost of each approximation step. To this end, for all quadrature points, we multiply the cost of solving the PDE with the cost of the truncated problem. 
As the truncation level of the random field is not fixed, we have to account for this change in the dimension. It is clearly more expensive to evaluate a diffusion coefficient which is truncated after a large number of terms than a diffusion coefficient involving only a few parameters. We therefore measure the cost of the truncated problem in terms of the active parameters, but neglect the possible cost of computing the KL. Similar to the cost models introduced in \cite{dick_high-dimensional_2013, kuo_liberating_2010} for infinite dimensional integration, two choices seem reasonable: Either the cost computing the truncated diffusion coefficient depends on the highest index of active variables, i.e. 
\begin{align} \label{eq:costA} \cost_{\text{Trunc}}=
	 \max \{k\in \mathbb{N}: l_{y_k}\geq 1\} \end{align}
or is measured by the number of active variables, i.e.
\begin{align}
 \cost_{\text{Trunc}}= \#\{k \in \mathbb{N}: l_{y_k}\geq 1\}.
\end{align}
The cost of the finite element method is measured by the number of finite element nodes. Here we assume for reasons of simplicity an optimal solver like multigrid or a multilevel preconditioner, such that only a cost of  $\mathcal{O}\left((N^{(\xx)}_{l_{\xx}})^n\right)$ is involved. And lastly, the number of quadrature points is used to capture the cost of the quadrature method, such that the total cost of computing $\bm{P}_{\tilde{\bm{l}}}u$ is given by \begin{align}\label{eq:cost}
{c_{\tilde{\bm{l}}}=(N^{(\xx)}_{l_{\xx}})^n \; \cdot \; \prod_{i=1}^m N^{(y_i)}_{l_{y_i}} \; \cdot \; \cost_{\text{Trunc}}(\bm{l}_{\yy}). }\end{align}

This way, we balance the error with the required work and our algorithm selects the index from $\tilde{\mathcal{A}}$ with the best benefit-cost ratio 
and minimizes $\sum_{\tilde{\bm{l}}\notin \tilde{\mathcal{I}}}  \|\Delta_{\tilde{\bm{l}}}u \|_\mathcal{X}$ for a given cost. Here,  the index with the highest contribution is only searched along the direct neighborhood of $\tilde{\mathcal{O}}$. Hence, we have to assume that the contributions of $\Delta_{\tilde{\bm{l}}} u$ decrease with the size of $\tilde{\bm{l}}$ to be able to ensure  convergence. This is comparable to the saturation assumption for adaptive finite element methods.
Having the coefficients in the Karhunen-Lo\`eve expansion sorted in decreasing order with respect to the $L^\infty(D)$-norm implies already  a decay for the truncation direction.  Note that this saturation assumption is not met, for example, when the diffusion field is rough  and contains discontinuities or spikes on a very fine scale. But in such a case other methods such as multilevel Monte Carlo or Quasi Monte Carlo would also run into similar problems.

In practice, this saturation assumption cannot be verified without additional knowledge. If exact knowledge on the contributions of the $\Delta_{\tilde{\bm{l}}} u$ would be available, an a-priori construction of an optimal index $\tilde{\mathcal{I}}$ could be used which would reduce the computation excess. But this information is often not available, e.g. for different random fields, and decay rates and corresponding optimal index sets have to be verified for each problem anew. Our a-posteriori algorithm however can be directly applied without any a-priori knowledge.
Note to this end that the initial spatial discretization has the be carefully chosen: An insufficiently-resolved mesh might not capture all essential features of the PDE problem and the algorithm may stop too early. On the other hand a too fine initial discretization would increase the overall cost. 
For our experiments we used an heuristic  to determine a good initial discretization. It is based on  an estimate/knowledge of the correlation length of the random field. More strategies can be found in \cite{giles_multilevel_2015}.

At last, let us mention different stopping criteria.
The criterion that is used in Algorithm \ref{alg:dimadaptKTinfinite} is based on a global profit indicator in the neighborhood set $\tilde{\mathcal{A}}$. 
This approach estimates how much profit is not covered yet. As the error contribution for indices
in the complement of $\tilde{\mathcal{I}}$ is not known, the only available information are the
contributions of the indices in $\tilde{\mathcal{A}}$.
Therefore, an estimate of the global profit is 
\begin{align*}
\eta=\sum_{\tilde{\bm{l}}\in \tilde{\mathcal{A}}} \eta_{\tilde{\bm{l}}}
\end{align*} and the algorithm stops when the global profit indicator is below a tolerance $\varepsilon$.
 It is also possible to replace  in the stopping criterion the profit by
the error contribution providing a global error indicator. 
Alternatively, the iteration could be stopped when a  given upper bound on the total work $C =\sum_{\tilde{\bm{l}}\in \tilde{\mathcal{I}}} c_{\tilde{\bm{l}}}$
is exhausted. Another stopping criterion would look at the difference between the solutions of two
iterations in order to estimate the error, i.e.,
$\| P_{\tilde{\mathcal{I}}^{(t+1)}}^{\text{ct}} u - P_{\tilde{\mathcal{I}}^{(t)}}^{\text{ct}} u \|_{\mathcal{X}} \leq \sum_{\tilde{\bm{l}}\in \tilde{\mathcal{I}}^{(t+1)}\setminus \tilde{\mathcal{I}}^{(t)}}
\|\Delta_{\tilde{\bm{l}}} u \|_\mathcal{X}.
$
The algorithm stops if this difference is $\zeta$-times in a row smaller than
$\varepsilon$.
Here, the parameter $\zeta$ is introduced to avoid stopping too early, when the solution differs
only slightly from one iteration to another while the next solution would result in a larger
difference.


\section{Numerical examples} \label{sec:numex}

As indicated in the introduction,  we focus on the model problem 
\begin{align*} 
	-\divergence (a(\xx,\omega) \nabla u(\xx, \omega))=f(\xx) \quad \text{ in } D,
\end{align*}
where $a(\xx,\omega)= \exp(b(\xx, \omega))$, compare  \eqref{eq:defalognormal}.
Moreover, we set $n=2$, $D=[0,1]^2$, and apply homogeneous Dirichlet boundary
conditions.
For our numerical examples, we assume that the covariance function of the random
field $b(\xx,\omega)$ is in the class of Mat\'ern covariance functions (see,
e.g. \cite{rasmussen_gaussian_2006}).
Then, the regularity and hence the eigenvalue decay for the Karhunen-Lo\`eve
expansion is  explicitly known. 

\begin{figure}
	\centering
	\includegraphics[width=0.48\textwidth]{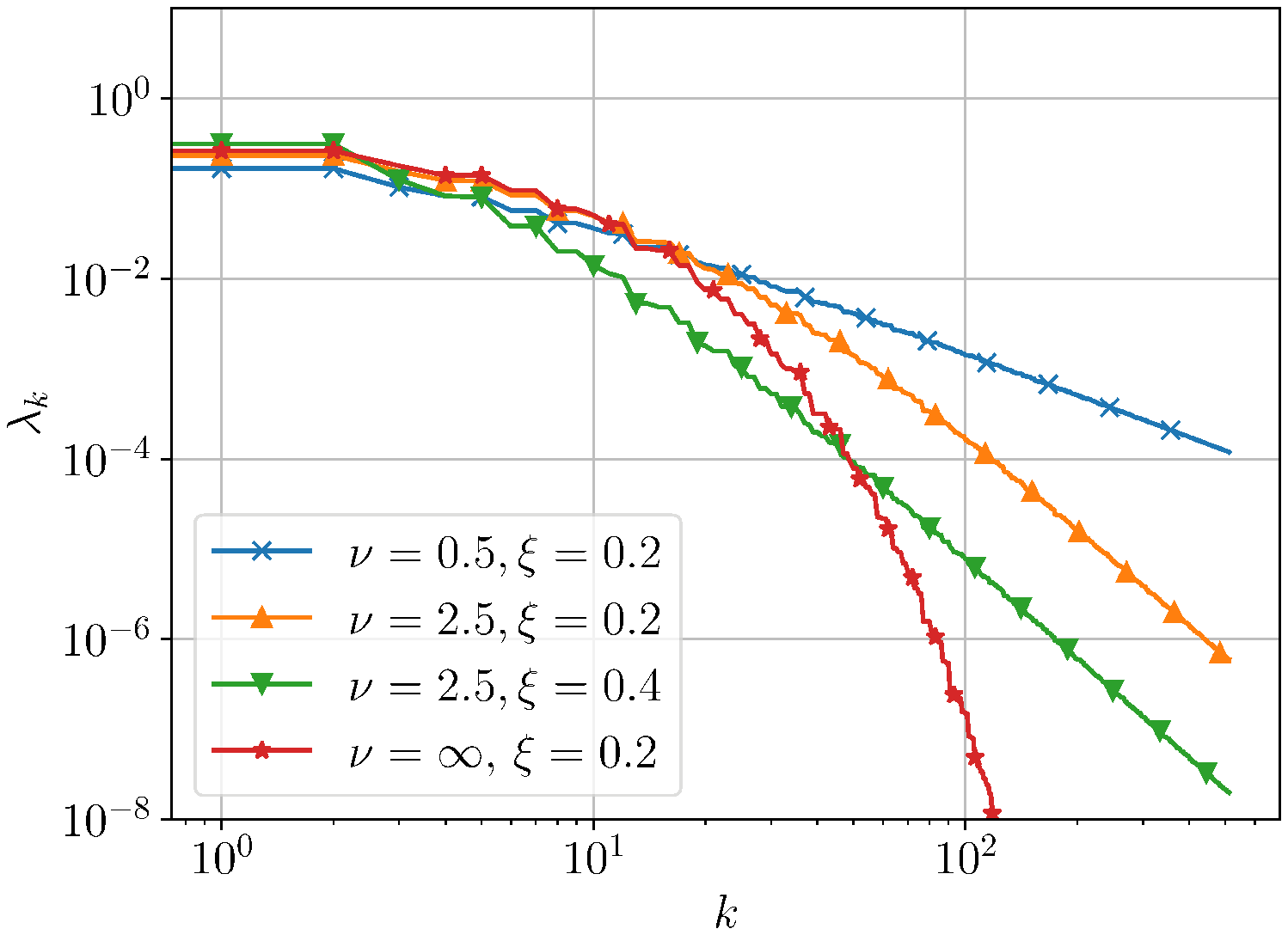}
	\includegraphics[width=0.48\textwidth, trim={0 0 1.5cm 0 },
	clip]{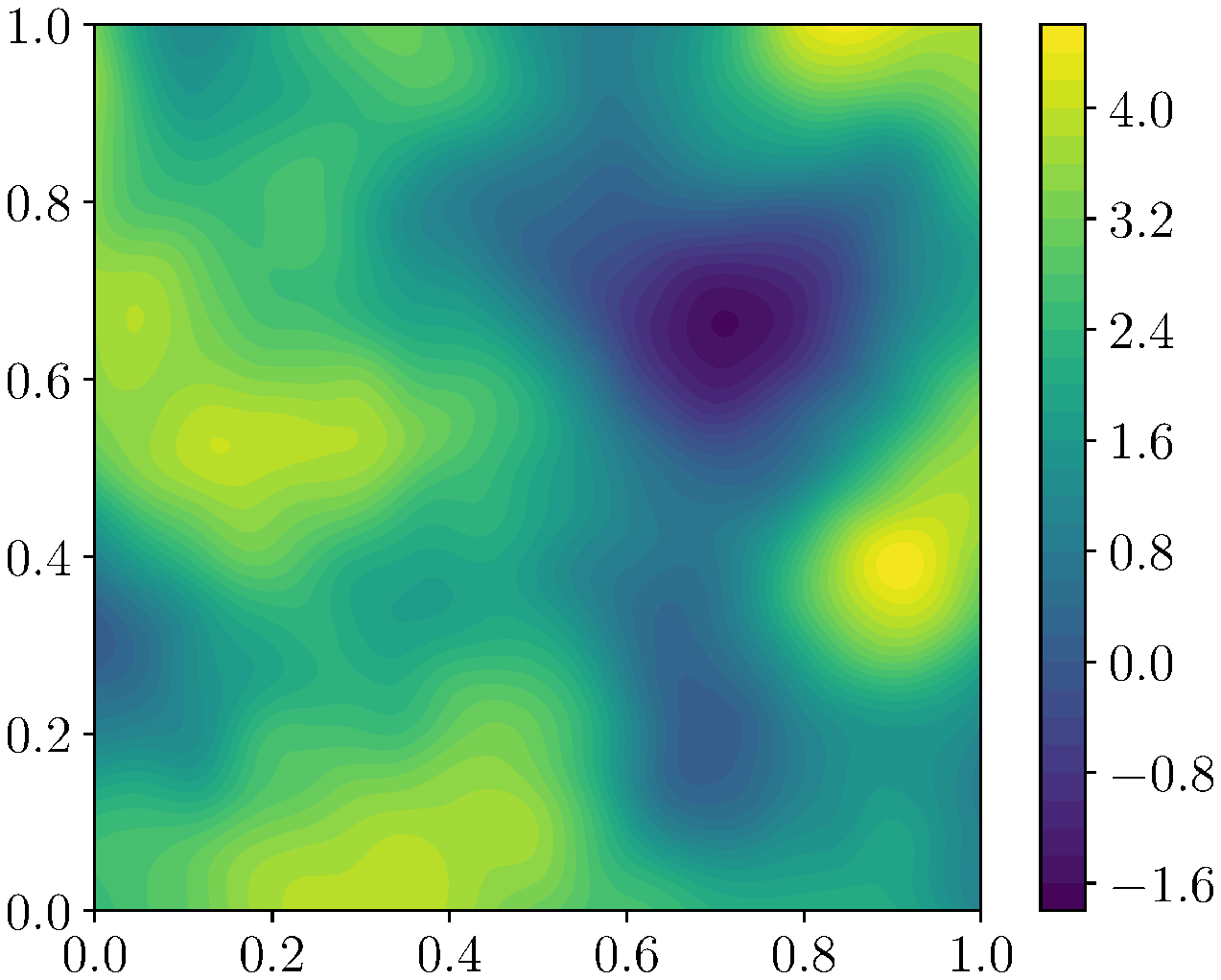}
	\caption{{Left: Decay of Karhunen-Lo\`eve coefficients $\lambda_k$ for
			different parameters $\nu$ and $\xi$. Right: Realization of random field
			$b(\xx,\omega)$ for $\xi=0.2$ and $\nu=2.5$.}}
	\label{fig:decaysequences}
\end{figure}

The Mat\'ern covariance kernels are defined by a stationary covariance function.
Let \linebreak $r=|\xx-\xx'|$ denote the distance between two points and  let
$\xi,\sigma^2>0$ be the correlation length and variance, respectively.  Then the
Mat\'ern covariance function  of order $\nu>0$ is defined as
\begin{align}\label{eq:materncov}
	\cov(r; \nu,\xi) =
	\sigma^2\frac{2^{1-\nu}}{\Gamma(\nu)}\Bigg(\sqrt{2\nu}\frac{r}{\xi}\Bigg)^\nu
	K_\nu\Bigg(\sqrt{2\nu}\frac{r}{\xi}\Bigg),
\end{align}
where $\Gamma$ is the Gamma function and $K_\nu$ is the modified Bessel function
of the second kind of order $\nu$.
The formula \eqref{eq:materncov} simplifies for $\nu=s+\frac{1}{2}$ with
$s\in\mathbb{N}$  to 
\begin{align*}
	\cov\left(r;s+\frac{1}{2},\xi\right) =
	\sigma^2\exp\left(-\frac{\sqrt{2s+1}r}{\xi}\right)\frac{s!}{(2s)!}\sum_{i=0}^s\frac{(s+i)!}{i!(s-i)!}\left(\frac{2\sqrt{2s+1}r}{\xi}\right)^{s-i}.
\end{align*}
Furthermore, in the limit $\nu \to \infty$, we obtain a Gaussian covariance
function 
\begin{align*}
	\cov(r; \xi )=\sigma^2\exp\left(-\frac{1}{2}\frac{r^2}{\xi^2}\right).
\end{align*}
The parameter $\nu$ can be regarded as a smoothness parameter as it  controls
the regularity of the covariance function at the point $r=0$ as well as the
regularity in the spatial variable of realizations of the random field. 
Both parameters $\xi$ and $\nu$ influence the decay of the Karhunen-Lo\`eve
coefficients and therefore the stochastic regularity. 
In the case $\nu < \infty $ the asymptotic decay rates of the KL coefficients
are known to be algebraic, while for $\nu=\infty$ the coefficients decrease
exponentially. 
Moreover, their pre-asymptotic behavior is determined by  the correlation length
$\xi$. The decay of the eigenvalues for various parameters are shown in Figure
\ref{fig:decaysequences} together with  a realization of the random field
$\log(a(x,\omega))$ for $\nu=2.5$ and $\xi=0.2$. Here, we obtained the KL
expansion by solving a discretized eigenvalue problem where the discretization
is much finer than the arising spatial discretizations for the PDE.

In the following we consider different modeling parameters $\nu$ and $\xi$
leading to different decay rates in the Karhunen-Lo\`eve expansion. Even though
the decay rates differ, we can apply the adaptive combination technique without
any adjustments, which would be necessary for any a-priori construction of a
sparse grid approximation.

We illustrate the performance of our new dimension-adaptive combination
technique by computing the first moment $\mathbb{E}[u]$ for different
permeability fields. For other quantities of interest such as the second moment
and linear functionals of $u$, we observed similar results in further
experiments. Here, we only present the results for the first moment in detail.
We assume homogeneous Dirichlet boundary condition and set $f=1$. 
The algorithm is applied with $\widehat{M}=5$ variables and the error is
measured in the $L^2(D)$-norm   with respect to a reference solution which was
obtained by using the combination technique with much finer discretizations
levels. 
Consequently, we use as a measure for the local contribution 
$\|\Delta_{\tilde{\bm{l}}} u\|_{L^2(D)}$ throughout the algorithm. The cost model is chosen as described in Section \ref{sec:profit} with \eqref{eq:costA} as this choice for the cost of the truncation better captures the actual computational time.
For solving the PDE for fixed parameter values on a uniform triangulation of $D$
with mesh width $h_{l_{\xx}}= \mathcal{O}(2^{-l_{\xx}})$ and piecewise finite
elements, the finite element software FEniCS \cite{alnaes_fenics_2015} is used. 
The experiments were done on a 3GHz Xeon 6136 CPU with 384GB of RAM.

\begin{figure}
	\centering
		\includegraphics[width=0.48\textwidth]{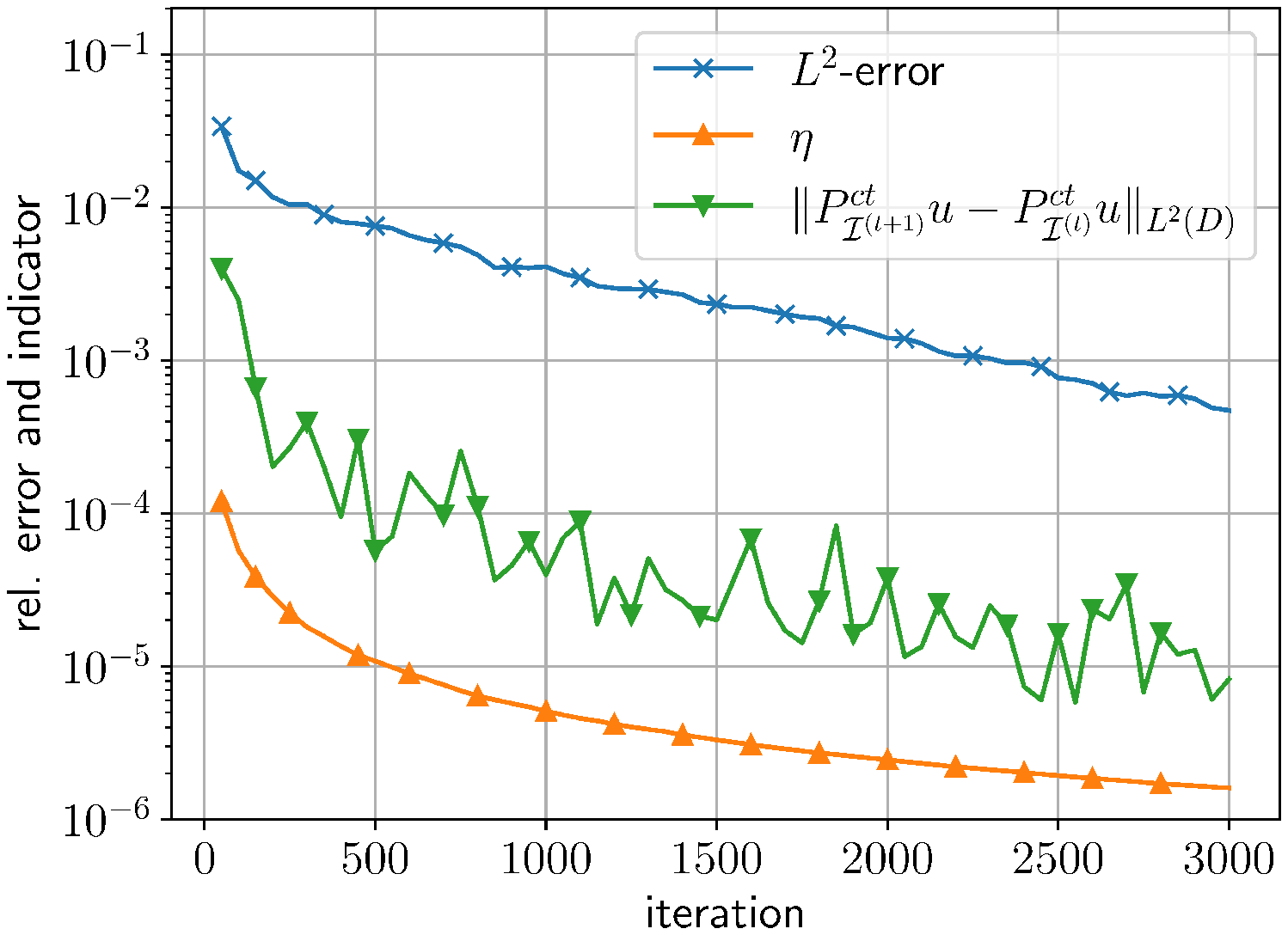}
		\includegraphics[width=0.48\textwidth]{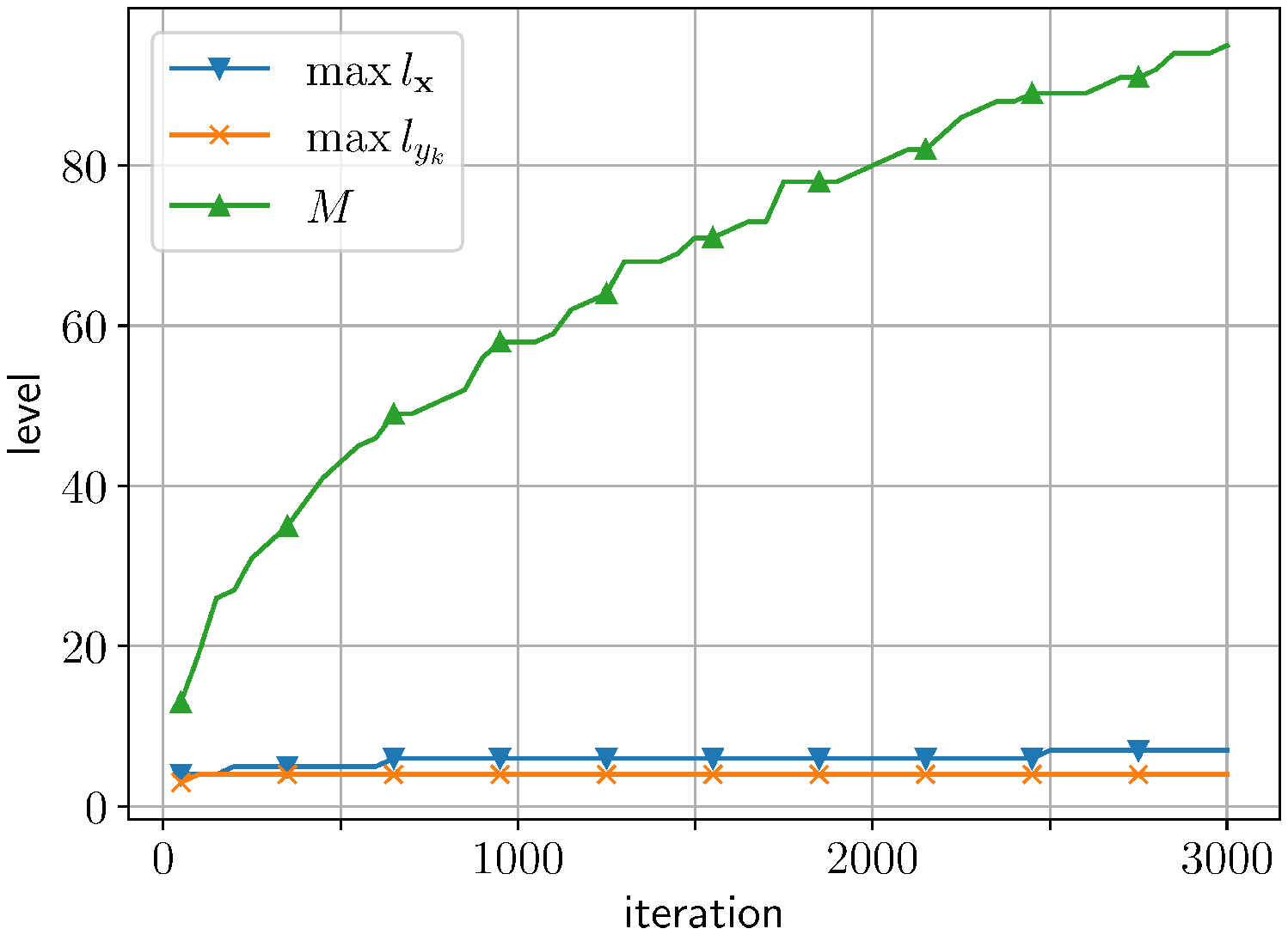}
	\caption{{Left: convergence in terms of number of iterations. Right: maximal
			level in the index set for spatial refinement, quadrature and truncation level.
	}}
	\label{fig:pb1detailed}
\end{figure}

The behavior of the dimension-adaptive algorithm for the Mat\'ern covariance
function with parameter $\nu=2.5$ and correlation length $\xi=0.4$ is
illustrated in Figure \ref{fig:pb1detailed}.
We observe that the error decreases  with increasing number of iterations. We
observe also that the error decreases fast in the first iterations as the
profits are high for these indices and that the global error estimators follow
this behavior of the error.

Figure \ref{fig:pb1detailed} (right) expresses the interplay of the three
different numerical approximations and the evolution of the index set during the
adaptive algorithm. It shows how the maximal level for each refinement direction
increases throughout the process, whereby, for the quadrature, the maximal level
of all dimensions is considered. As expected, a finer discretization for the
spatial domain than for the stochastic parameters  is needed. For example, after
$1000$ iterations, the highest included  spatial level is $l_{\xx}=6$ which
corresponds to $257$ degrees of freedom in one direction for the finite element
function. In contrast, at most $16$ quadrature points are used for the
univariate quadrature rules, since the Gaussian quadrature exhibits exponential
convergence and therefore not many levels are required for a small quadrature
error.
Moreover, the number of activated parameter dimensions increases successively 
such that, after 3000 iterations, our algorithm has activated $95$ parameter
dimensions which cover the variability of the random field to a large extent.
Moreover, it balances the corresponding finite noise truncation error with the
quadrature error  and the spatial discretization error while also taking the
associated costs into account.

\begin{figure}
	\centering
	\includegraphics[width=0.5\textwidth, trim={3cm 0cm 0 1cm},
	clip]{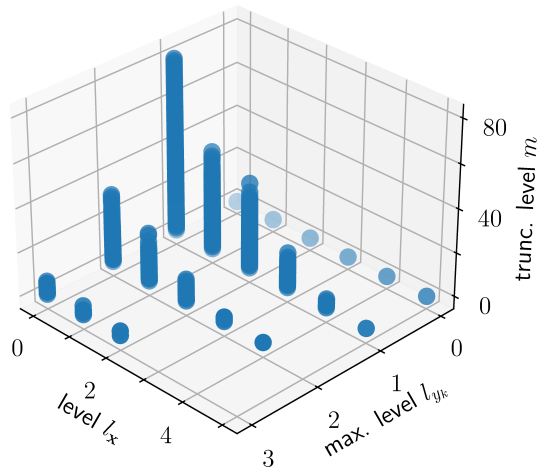}
	\caption{{Index set $\tilde{\mathcal{O}}$ at iteration $2000$.}}
	\label{fig:it2000}
\end{figure}

The interplay  of the truncation, the quadrature and the finite element method
is also illustrated in Figure \ref{fig:it2000}, where the index set
$\tilde{\mathcal{O}}$ at iteration $2000$ is depicted. Again, we
give  just the maximal level over all quadrature rules for the
stochastic approximation direction.  We observe the typical shape  of a
generalized sparse grid discretization in three dimensions: A high resolution in
one direction relates to a corresponding  low resolution in the other two
directions. 

\begin{figure}
	\centering
	\includegraphics[width=0.48\textwidth,trim={0 0 0 -0.6cm},
	clip]{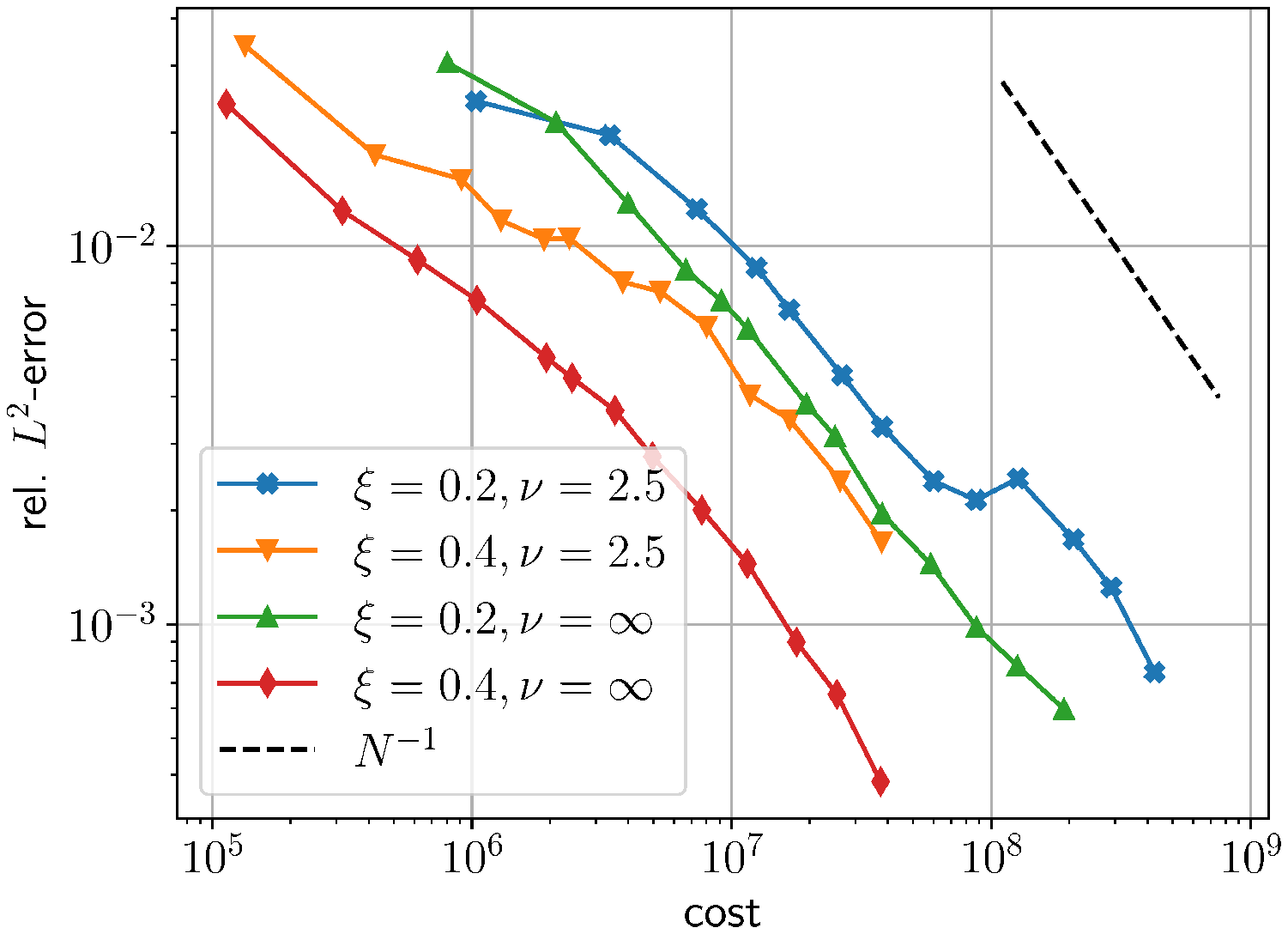}
	\includegraphics[width=0.48\textwidth,trim={0 0 0 -0.6cm},
	clip]{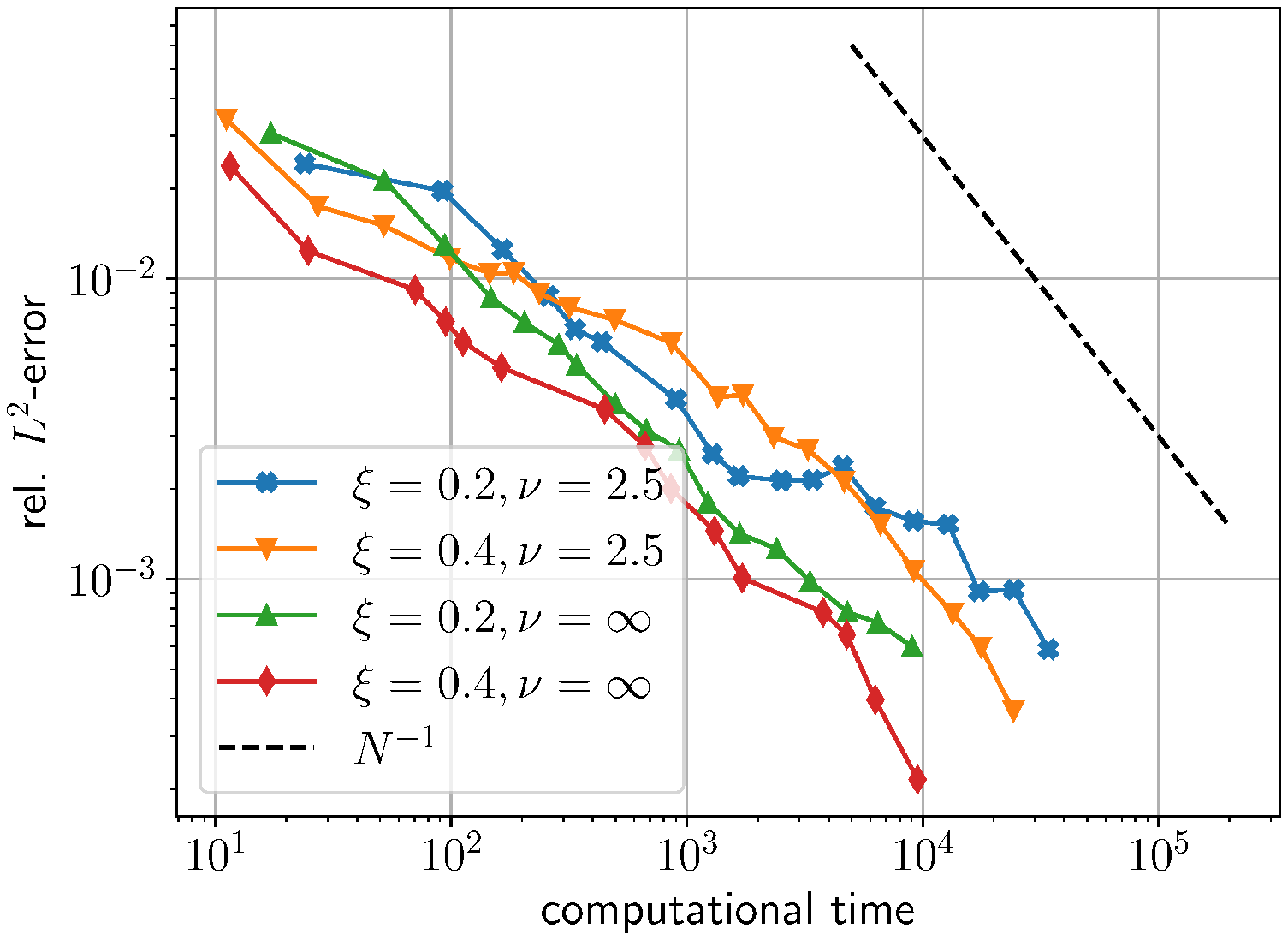}
	\caption{{Convergence of dimensions-adaptive {algorithm} with respect to the
			computational cost (left) and computational time (right).}}
	\label{fig:convergence}
\end{figure}

The convergence rate of the solution with respect to the computational cost
 and computational time is shown in Figure \ref{fig:convergence} for the
example of random fields with the parameter values $\nu=2.5$ and $\nu=\infty$
and the correlation lengths $\xi=0.2$ and $\xi=0.4$ which possess an algebraic
and an exponential decay of the coefficients, respectively. In all four cases,
the decay of the coefficients is sufficiently fast such that the mixed
regularity condition as in \cite{harbrecht_multilevel_2016} for a sparse grid
approximation between the spatial and the parametric variables is met.
We can thus expect that the convergence rate of our sparse grid approach
is limited by the slowest approximation direction as mentioned in Section \ref{sec:SG}. 
As the realizations of the random fields are continuously differentiable,  the
PDE solution satisfies  $u\in L^q(\Omega, H^2(D))$ for $q<\infty$ (see
\cite{charrier_finite_2013,teckentrup_further_2013}). Hence, for fixed $\yy$,
the finite element method achieves at most a convergence rate of order one with
respect to the number of degrees of freedom (cf. Section \ref{sec:fem}), as we
have $\mathcal{O}((N^{(\xx)}_{l_{\xx}})^2)$ degrees of freedom for a mesh with
width $h_{l_{\xx}}=(N^{(\xx)}_{l_{\xx}})^{-1}$. The truncation error in terms of
$m$ decreases with a rate faster than one and the smooth dependence on the
stochastic parameters yields also a faster convergence of the quadrature. 
Hence the convergence rate for smooth fields is limited by that of the  finite
element method. 
Indeed,  we observe an asymptotic convergence rate that approaches
$\mathcal{O}(N^{-1})$ with respect to {the overall computational cost measured
	by $N= \sum_{\tilde{\bm{l}}\in\tilde{\mathcal{I}}} c_{\tilde{\bm{l}}}$} (c.f.
\eqref{eq:cost}), as well as with respect to the computational time.

This is consistent with the behavior expected from a-priori analysis and
other adaptive algorithms. For example, the results in
\cite{haji-ali_multi-index_2016} also show that the spatial discretization
dominates the convergence for smooth truncated random field. 
Note that, compared to the algorithm \cite{owen_dimension-adaptive_2018} which
controls the spatial and truncation level in a multi-index Monte Carlo method,
our method achieves a higher convergence rate by a factor of about two. Since our approach with the
Gaussian quadrature can exploit the regularity  of the parametric variables, the
convergence rate is not limited to $1/2$ as for Monte Carlo methods and now the
spatial regularity  becomes the limiting factor. 
	
\begin{figure}
	\centering
	\includegraphics[width=0.48\textwidth,trim={0 0 0 -0.6cm},
	clip]{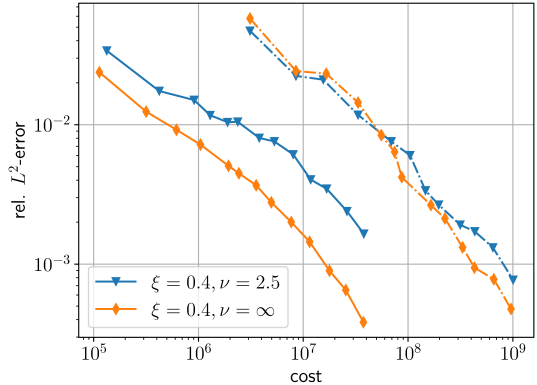} 
	\caption{{Convergence of our dimensions-adaptive {algorithm} (solid line) and
			of the algorithm with fixed spatial discretization (dashed line ) with respect
			to the computational cost.}}
	\label{fig:fixedlx}
\end{figure}

To illustrate the benefit of the combination of all three discretization steps,
we compare our results to those in \cite{garcke_adaptive_2016, chen_sparse_2018, ernst_convergence_2018} where the spatial discretization is fixed to a very fine
level. To this end, we solve the PDE problems on a mesh with $257\times 257$ nodes
which is sufficiently refined. This indeed is the approach which was followed in \cite{garcke_adaptive_2016}, albeit with a different resolution and for a different quantity of interest.  Figure  \ref{fig:fixedlx} shows the convergence
rate for two random fields with respect to the computational cost which also
includes the cost of the spatial discretization. As it can be expected, using a
very fine spatial discretization involves much higher cost for the same error tolerance as our algorithm because the discretization is not balanced with the parametric approximation.

\begin{figure}
	\centering
	\includegraphics[width=0.88\linewidth, trim={0 0.25cm 0 0}, clip
	]{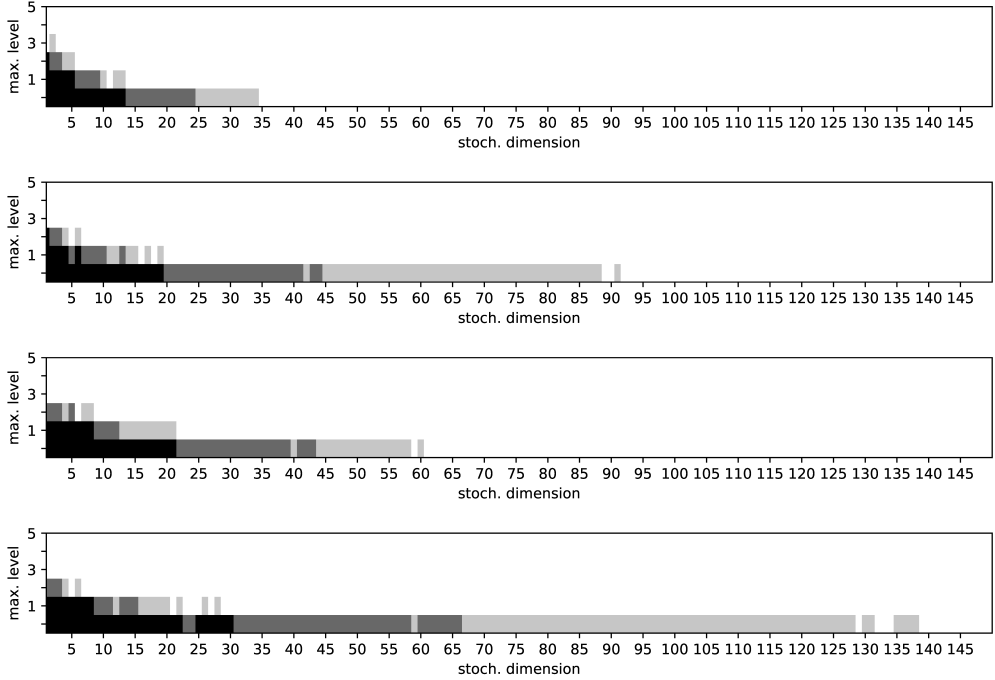}
	\caption{{Maximal quadrature level (black after iteration $100$, dark gray
			after iteration $500$ and light gray after iteration $2500$) for different
			random fields. From top to bottom: $\nu=\infty$ and $\xi=0.4$; $\nu=2.5$ and
			$\xi=0.4$; $\nu=\infty$ and $\xi=0.2$; $\nu=2.5$ and  $\xi=0.2$.}}
	\label{fig:parametericdisc}
\end{figure}

In Figure \ref{fig:parametericdisc} we plot the levels used for the
parametric discretization to illustrate the adjustment of our algorithm to the
specific problem under consideration.
Here, the largest level for each one-dimensional quadrature level is shown. We
observe that a higher refinement is applied for the first few parameter
dimensions as their influence on the solution is larger, while for the later
parameter dimensions only the lowest level quadrature is required. Our algorithm
therefore detects the anisotropy in the parameters and adjusts the necessary
quadrature levels accordingly. 
We also observe that fewer  variables are activated for the smoother  fields
with $\nu=\infty$ than for the fields with $\nu=2.5$. 
Due to the faster decay of the KL coefficients, it suffices to only consider a
smaller number of parametric variables. 
The correlation length merely affects the size of the pre-asympototic regime of
the decay. For smaller correlation lengths the regime is larger, which is
reflected in the quadrature levels by a larger number of parametric dimensions
that are refined to a higher extent. Furthermore, more variables are
activated in the computation. However, the asymptotic convergence rate is not
affected. 

As a final example, we consider a problem with $\nu=0.5$ and thus with an
exponential covariance function. 
For $r=|\xx-\xx'|_1$ the Karhunen-Lo\`eve eigenvalues and eigenfunctions are
explicitly known. 
In this case the almost sure convergence of $a(\xx, \omega)$ in $L^\infty(D)$ is
guaranteed \cite{charrier_strong_2012}, but the condition $\{\gamma_k\}_k
\in\ell^1(\mathbb{N})$ (c.f. \eqref{eq:defgamma}) is not satisfied as the
realizations of the diffusion coefficient are only H\"older continuous with
exponent $1/2$. 
For this example, the assumptions typically made to obtain a-priori estimates on
the size of the surpluses $\Delta_{\tilde{\bm{l}}} u$ cannot be verified.
Nevertheless, the  set that includes the indices with the highest local
benefit-cost ratio can be constructed with our adaptive algorithm in a
straightforward way. 

\begin{figure}
	\centering
	\includegraphics[width=0.48\linewidth]{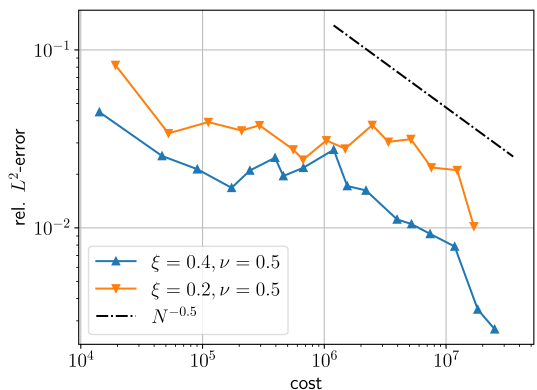}
	\caption{Convergence plot for $\nu=0.5$}
	\label{fig:convergence0.5}
\end{figure}

The rougher diffusion coefficient results in  less regular solutions of the PDE
which only  satisfies \linebreak $u(\yy)\in H^{3/2}(D)$, see
\cite{charrier_finite_2013,teckentrup_further_2013}. 
Hence, we can only  expect from the finite element method to achieve a
convergence order of $1/2$. Indeed, we observe this slower decay of the
$L^2(D)$-error with our  dimension-adaptive combination technique in Figure
\ref{fig:convergence0.5}. 

The associated slow decrease of the coefficients in the KL expansion results in
a large number of parametric variables that need to be taken into account. Our
algorithm  adjusts to this slower decay and activates $192$ parametric variables
for $\xi=0.4$ at iteration $2500$ (and $189$ parametric variables for $\xi=0.2$
at iteration $2500$) which  cover most of the variability of the random field.
As the parametric dimension grows and more parameters have a strong influence on
the solution, the sparse grid for the parametric discretization gets naturally
more costly than for problems with only a few important parametric dimensions.

\section{Concluding remarks}
In this paper we presented a dimension-adaptive algorithm that computes quantities of interests for a given stochastic random diffusion coefficient and relies on a sparse grid approximation. The algorithm balances the spatial, stochastic and parametric approximation  in a cost-efficient way and detects the necessary discrectization levels to obtain a prescribed overall error.
This method is especially useful when information about the decrease of the increments $\Delta_{\tilde{\bm{l}}}u$  is not available, as it does not require any a-priori knowledge on the regularity of the solution and on the decay of the KL coefficients. 

Our adaptive algorithm is based on the benefit-cost ratio as reliable refinement indicator such that it minimizes the error for a prescribed computational cost. Further work may focus on the optimality of this refinement strategy and on the choice of efficient and reliable global error estimators. 

Numerical examples showed that our  obtained convergence order is comparable with the known a-priori convergence results for smooth random diffusion coefficients. Compared to other methods such as multilevel Monte Carlo and Quasi Monte Carlo, our algorithm can exploit the parametric regularity such that the convergence order in terms of  the computational cost is limited by the finite element method for smooth fields. For problems with less smoothness, the algorithm also yields good results, but converges slower due to the lower regularity of the PDE solution.

As the algorithm is based on combinations of quadratures and finite element discretizations, the algorithm could be generalized using other spatial discretization techniques, such as h-adaptive and h-p-adaptive finite element methods. By this, we
expect to be able to benefit from the regularity in the case of smooth PDE solutions and to obtain a higher order of convergence than by a Quasi Monte Carlo method.  
Furthermore, our algorithm is not limited to computing first and second moment or linear functionals. In fact, it can easily be  generalized to other quantities of interest involving a different function $\mathcal{F}$. This might however require an additional approximation step for the evaluation of $\mathcal{F}$ which also needs to be balanced. 

In this work, we concentrated on adaptivity with respect to the stochastic parameters and kept the spatial discretization uniform. Besides, there is also the possibility/necessity to apply local adaptive refinement in space, compare e.g. \cite{bespalov_error_2021, bespalov_error_2022,feischl_convergence_2021}.
The study of adaptivity in the stochastic {\em and} the spatial domain is future work.


\section*{Acknowledgments}
This research was supported by the\textit{ Hausdorff Center for Mathematics} in Bonn, funded by the Deutsche Forschungsgemeinschaft (DFG, German Research Foundation) under Germany’s Excellence Strategy - EXC-2047/1 - 390685813 and the CRC 1060 \textit{ The Mathematics of Emergent Effects} of the Deutsche Forschungsgemeinschaft.











\bibliographystyle{siamplain}
\bibliography{Lit}

\begin{thebibliography}{10}

\bibitem{ainsworth_posteriori_2000}
{\sc M.~Ainsworth and J.~Oden}, {\em A Posteriori Error Estimation in Finite
  Element Analysis}, Pure and Applied Mathematics, Wiley, New York, 2000.

\bibitem{alnaes_fenics_2015}
{\sc M.~Alnaes, J.~Blechta, J.~Hake, A.~Johansson, B.~Kehlet, A.~Logg,
  C.~Richardson, J.~Ring, M.~Rognes, and G.~Wells}, {\em The {FEniCS} {Project}
  {Version} 1.5}, Archive of Numerical Software, Vol 3 (2015).

\bibitem{babuska_stochastic_2010}
{\sc I.~Babuška, F.~Nobile, and R.~Tempone}, {\em A {stochastic} {collocation}
  {method} for {elliptic} {partial} {differential} {equations} with {random}
  {input} {data}}, SIAM Rev., 52 (2010), pp.~317--355.

\bibitem{bachmayr_sparse_2017}
{\sc M.~Bachmayr, A.~Cohen, R.~DeVore, and G.~Migliorati}, {\em {Sparse
  polynomial approximation of parametric elliptic {PDEs}. {Part} {II}:
  Lognormal coefficients}}, ESAIM Math. Model. Numer. Anal., 51 (2017),
  pp.~341--363.

\bibitem{barth_multi-level_2011}
{\sc A.~Barth, C.~Schwab, and N.~Zollinger}, {\em {Multi-level {Monte} {Carlo}
  {finite} {element} method for elliptic {PDEs} with stochastic coefficients}},
  Numer. Math., 119 (2011), pp.~123--161.

\bibitem{bear_modeling_2010}
{\sc J.~Bear and A.~Cheng}, {\em Modeling {Groundwater} {Flow} and
  {Contaminant} {Transport}}, Springer Netherlands, 2010.

\bibitem{beck_optimal_2012}
{\sc J.~Beck, R.~Tempone, F.~Nobile, and L.~Tamellini}, {\em On the optimal
  polynomial approximation of stochastic {PDEs} by {Galerkin} and collocation
  methods}, Mathematical Models and Methods in Applied Sciences, 22 (2012),
  \url{https://doi.org/10.1142/S0218202512500236}.

\bibitem{bespalov_error_2022}
{\sc A.~Bespalov and D.~Silvester}, {\em Error estimation and adaptivity for
  stochastic collocation finite elements, {Part} {II}: Multilevel
  approximation}, SIAM Journal on Scientific Computing, 45 (2023),
  pp.~A781--A797, \url{https://doi.org/10.1137/22M1479361}.

\bibitem{bespalov_error_2021}
{\sc A.~Bespalov, D.~J. Silvester, and F.~Xu}, {\em Error estimation and
  adaptivity for stochastic collocation finite elements, {Part} {I}:
  Single-level approximation}, SIAM Journal on Scientific Computing, 44 (2022),
  pp.~A3393--A3412, \url{https://doi.org/10.1137/21M1446745}.

\bibitem{bungartz_sparse_2004}
{\sc H.-J. Bungartz and M.~Griebel}, {\em Sparse {grids}}, Acta Numer., 13
  (2004), pp.~147--269.

\bibitem{charrier_strong_2012}
{\sc J.~Charrier}, {\em Strong and {weak} {error} {estimates} for {elliptic}
  {partial} {differential} {equations} with {random} {coefficients}}, SIAM J.
  Numer. Anal., 50 (2012), pp.~216--246.

\bibitem{charrier_weak_2013}
{\sc J.~Charrier and A.~Debussche}, {\em Weak truncation error estimates for
  elliptic {PDEs} with lognormal coefficients}, Stoch. Partial Differ. Equ.
  Anal. Comput., 1 (2013), pp.~63--93.

\bibitem{charrier_finite_2013}
{\sc J.~Charrier, R.~Scheichl, and A.~Teckentrup}, {\em Finite {element}
  {error} {analysis} of {elliptic} {PDEs} with {random} {coefficients} and
  {its} {application} to {multilevel} {Monte} {Carlo} {methods}}, SIAM J.
  Numer. Anal., 51 (2013), pp.~322--352.

\bibitem{chen_sparse_2018}
{\sc P.~Chen}, {\em Sparse quadrature for high-dimensional integration with
  {Gaussian} measure}, ESAIM Math. Model. Numer. Anal., 52 (2018),
  pp.~631--657.

\bibitem{chkifa_high-dimensional_2014}
{\sc A.~Chkifa, A.~Cohen, and C.~Schwab}, {\em High-{dimensional} {adaptive}
  {sparse} {polynomial} {interpolation} and {applications} to {parametric}
  {PDEs}}, Found. Comput. Math., 14 (2014), pp.~601--633.

\bibitem{dick_high-dimensional_2013}
{\sc J.~Dick, F.~Y. Kuo, and I.~H. Sloan}, {\em High-dimensional integration:
  {The} quasi-{Monte} {Carlo} way}, Acta Numerica, 22 (2013), pp.~133--288,
  \url{https://doi.org/10.1017/S0962492913000044}.

\bibitem{dung_hyperbolic_2016}
{\sc D.~Dũng and M.~Griebel}, {\em Hyperbolic cross approximation in infinite
  dimensions}, Journal of Complexity, 33 (2016), pp.~55--88,
  \url{https://doi.org/10.1016/j.jco.2015.09.006}.

\bibitem{dung__2018}
{\sc D.~Dũng, M.~Griebel, V.~N. Huy, and C.~Rieger}, {\em
  {$\varepsilon$}-dimension in infinite dimensional hyperbolic cross
  approximation and application to parametric elliptic {PDE}s}, Journal of
  Complexity, 46 (2018), pp.~66--89,
  \url{https://doi.org/10.1016/j.jco.2017.12.001}.

\bibitem{ernst_convergence_2018}
{\sc O.~Ernst, B.~Sprungk, and L.~Tamellini}, {\em Convergence of {sparse}
  {collocation} for {functions} of {countably} {many} {Gaussian} {random}
  {variables} (with {application} to {elliptic} {PDEs})}, SIAM J. Numer. Anal.,
  56 (2018), pp.~877--905.

\bibitem{feischl_convergence_2021}
{\sc M.~Feischl and A.~Scaglioni}, {\em Convergence of adaptive stochastic
  collocation with finite elements}, Comput. Math. Appl., 98 (2021),
  pp.~139--156.

\bibitem{garcke_sparse_2012}
{\sc J.~Garcke}, {\em Sparse {grids} in a {nutshell}}, in Sparse {Grids} and
  {Applications}, J.~Garcke and M.~Griebel, eds., vol.~88 of Lect. Notes
  Comput. Sci. Eng., Springer Berlin Heidelberg, 2012, pp.~57--80.

\bibitem{gerstner_dimension_2003}
{\sc T.~Gerstner and M.~Griebel}, {\em Dimension {adaptive} {tensor} {product}
  {quadrature}}, Computing, 71 (2003), pp.~65--87.

\bibitem{giles_multilevel_2015}
{\sc M.~B. Giles}, {\em Multilevel {Monte} {Carlo} methods}, Acta Numerica, 24
  (2015), pp.~259--328, \url{https://doi.org/10.1017/S096249291500001X}.

\bibitem{graham_quasi-monte_2015}
{\sc I.~Graham, F.~Kuo, J.~Nichols, R.~Scheichl, C.~Schwab, and I.~Sloan}, {\em
  {Quasi-{Monte} {Carlo} finite element methods for elliptic {PDEs} with
  lognormal random coefficients}}, Numer. Math., 131 (2015), pp.~329--368.

\bibitem{griebel_construction_2013}
{\sc M.~Griebel and H.~Harbrecht}, {\em On the construction of sparse tensor
  product spaces}, Mathematics of Computation, 82 (2013), pp.~975--994,
  \url{https://doi.org/10.1090/S0025-5718-2012-02638-X}.

\bibitem{griebel_multilevel_2020}
{\sc M.~Griebel, H.~Harbrecht, and M.~Multerer}, {\em Multilevel {quadrature}
  for {elliptic} {parametric} {partial} {differential} {equations} in {case} of
  {polygonal} {approximations} of {curved} {domains}}, SIAM J. Numer. Anal., 58
  (2020), pp.~684--705.

\bibitem{Griebel.Schneider.Zenger:1992}
{\sc M.~Griebel, M.~Schneider, and C.~Zenger}, {\em {A combination technique
  for the solution of sparse grid problems}}, in {Iterative Methods in Linear
  Algebra}, P.~de~Groen and R.~Beauwens, eds., IMACS, Elsevier, 1992,
  pp.~263--281.

\bibitem{haji-ali_novel_2018}
{\sc A.-L. Haji-Ali, H.~Harbrecht, M.~Peters, and M.~Siebenmorgen}, {\em {Novel
  results for the anisotropic sparse grid quadrature}}, J. Complexity, 47
  (2018), pp.~62--85.

\bibitem{haji-ali_multi-index_2016}
{\sc A.-L. Haji-Ali, F.~Nobile, L.~Tamellini, and R.~Tempone}, {\em
  Multi-{index} {stochastic} {collocation} for random {PDEs}}, Comput. Methods
  Appl. Mech. Engrg. C, 306 (2016), pp.~95--122.

\bibitem{harbrecht_multilevel_2016}
{\sc H.~Harbrecht, M.~Peters, and M.~Siebenmorgen}, {\em Multilevel
  {accelerated} {quadrature} for {PDEs} with {log}-{normally} {distributed}
  {diffusion} {coefficient}}, SIAM/ASA J. Uncertain. Quantif., 4 (2016),
  pp.~520--551.

\bibitem{hoang_n-term_2014}
{\sc V.~Hoang and C.~Schwab}, {\em N-{term} {Wiener} {chaos} {approximation}
  {rates} {for} {elliptic} {PDEs} {with} {lognormal} {Gaussian} {random}
  {inputs}}, Math. Models Methods Appl. Sci., 24 (2014), pp.~797--826.

\bibitem{jakeman_adaptive_2020}
{\sc J.~Jakeman, M.~Eldred, G.~Geraci, and A.~Gorodetsky}, {\em Adaptive
  multi‐index collocation for uncertainty quantification and sensitivity
  analysis}, International Journal for Numerical Methods in Engineering, 121
  (2020), pp.~1314--1343.

\bibitem{kuo_multilevel_2017}
{\sc F.~Kuo, R.~Scheichl, C.~Schwab, I.~Sloan, and E.~Ullmann}, {\em Multilevel
  quasi-{Monte} {Carlo} methods for lognormal diffusion problems}, Math. Comp.,
  86 (2017), pp.~2827--2860.

\bibitem{kuo_quasi-monte_2012}
{\sc F.~Kuo, C.~Schwab, and I.~Sloan}, {\em Quasi-{Monte} {Carlo} {finite}
  {element} {methods} for a {class} of {elliptic} {partial} {differential}
  {equations} with {random} {coefficients}}, SIAM J. Numer. Anal., 50 (2012),
  pp.~3351--3374.

\bibitem{kuo_liberating_2010}
{\sc F.~Y. Kuo, I.~H. Sloan, G.~W. Wasilkowski, and H.~Woźniakowski}, {\em
  Liberating the dimension}, Journal of Complexity, 26 (2010), pp.~422--454,
  \url{https://doi.org/10.1016/j.jco.2009.12.003}.

\bibitem{narayan_adaptive_2014}
{\sc A.~Narayan and J.~Jakeman}, {\em Adaptive {Leja} {sparse} {grid}
  {constructions} for {stochastic} {collocation} and {high}-{dimensional}
  {approximation}}, SIAM Journal on Scientific Computing, 36 (2014),
  pp.~A2952--A2983.

\bibitem{nobile_convergence_2016}
{\sc F.~Nobile, L.~Tamellini, and R.~Tempone}, {\em Convergence of
  quasi-optimal sparse-grid approximation of {Hilbert}-space-valued functions:
  Application to random elliptic {PDEs}}, Numerische Mathematik, 134 (2016),
  pp.~343--388, \url{https://doi.org/10.1007/s00211-015-0773-y}.

\bibitem{garcke_adaptive_2016}
{\sc F.~Nobile, L.~Tamellini, F.~Tesei, and R.~Tempone}, {\em An {adaptive}
  {sparse} {grid} {algorithm} for {elliptic} {PDEs} with {lognormal}
  {diffusion} {coefficient}}, in Sparse {Grids} and {Applications} -
  {Stuttgart} 2014, J.~Garcke and D.~Pflüger, eds., vol.~109, Springer
  International Publishing, 2016, pp.~191--220.

\bibitem{nobile_anisotropic_2008}
{\sc F.~Nobile, R.~Tempone, and C.~Webster}, {\em An {anisotropic} {sparse}
  {grid} {stochastic} {collocation} {method} for {partial} {differential}
  {equations} with {random} {input} {data}}, SIAM J. Numer. Anal., 46 (2008),
  pp.~2411--2442.

\bibitem{rasmussen_gaussian_2006}
{\sc C.~Rasmussen and C.~Williams}, {\em {Gaussian Processes for Machine
  Learning}}, Adaptive Computation and Machine Learning, MIT Press, 2006.

\bibitem{robbe_multi-index_2017}
{\sc P.~Robbe, D.~Nuyens, and S.~Vandewalle}, {\em A {Multi}-{index}
  {quasi}--{Monte} {Carlo} {algorithm} for {lognormal} {diffusion} {problems}},
  SIAM J. Sci. Comput., 39 (2017), pp.~S851--S872.

\bibitem{owen_dimension-adaptive_2018}
{\sc P.~Robbe, D.~Nuyens, and S.~Vandewalle}, {\em A dimension-adaptive
  multi-index {Monte} {Carlo} method applied to a model of a heat exchanger},
  in Monte {Carlo} and {Quasi}-{Monte} {Carlo} {Methods}, A.~B. Owen and P.~W.
  Glynn, eds., vol.~241, Springer International Publishing, Cham, 2018,
  pp.~429--445, \url{https://doi.org/10.1007/978-3-319-91436-7_24}.

\bibitem{schwab_sparse_2011}
{\sc C.~Schwab and C.~Gittelson}, {\em Sparse tensor discretizations of
  high-dimensional parametric and stochastic {PDEs}}, Acta Numer., 20 (2011),
  pp.~291--467.

\bibitem{teckentrup_multilevel_2015}
{\sc A.~Teckentrup, P.~Jantsch, C.~Webster, and M.~Gunzburger}, {\em A
  {multilevel} {stochastic} {collocation} {method} for {partial} {differential}
  {equations} with {random} {input} {data}}, SSIAM/ASA J. Uncertain. Quantif.,
  3 (2015), pp.~1046--1074.

\bibitem{teckentrup_further_2013}
{\sc A.~Teckentrup, R.~Scheichl, M.~Giles, and E.~Ullmann}, {\em Further
  analysis of multilevel {Monte} {Carlo} methods for elliptic {PDEs} with
  random coefficients}, Numer. Math., 125 (2013), pp.~569--600.

\end{thebibliography}
\end{document}